\documentclass[a4paper,11pt, reqno]{amsart}

\usepackage[DIV=12, oneside]{typearea} 
\usepackage[utf8]{inputenc}
\usepackage[T1]{fontenc}
\usepackage[english]{babel}
\usepackage{esint}
\usepackage{amssymb, amsthm, bbm, amsmath, latexsym, color}
\usepackage{enumerate, mathrsfs, hyphenat, graphicx}
\usepackage[font=sf, labelfont=sf]{caption}
\usepackage{thmtools, thm-restate}
\usepackage{thm-patch, thm-kv, thm-autoref}
\usepackage[colorlinks=true, urlcolor=blue, linkcolor=blue, citecolor=black]{hyperref}
\theoremstyle{plain}
\declaretheorem[title=Theorem, parent=section]{theorem}
\declaretheorem[title=Lemma,sibling=theorem]{lemma}
\declaretheorem[title=Corollary,sibling=theorem]{corollary}
\declaretheorem[title=Proposition,sibling=theorem]{proposition}
\theoremstyle{definition}
\declaretheorem[title=Definition,sibling=theorem]{definition}
\declaretheorem[unnumbered,title=Remark]{remark}

\numberwithin{equation}{section}
\newcommand{\eps}{\varepsilon}
\newcommand{\R}{  \mathbb{R}}

\newcommand{\N}{  \mathbb{N}}

\newcommand{\cB}{  \mathcal{B}}
\newcommand{\cC}{  \mathcal{C}}
\newcommand{\cD}{  \mathcal{D}}

\newcommand{\cF}{  \mathcal{F}}

\newcommand{\cL}{  \mathcal{L}}

\renewcommand{\epsilon}{\varepsilon}
\renewcommand{\d}{\, \textnormal{d}}

\addto\extrasenglish{%

}

\usepackage[backgroundcolor=none, bordercolor=blue]{todonotes}
\makeatletter
\providecommand\@dotsep{5}
\renewcommand{\listoftodos}[1][\@todonotes@todolistname]{%
  \@starttoc{tdo}{#1}}
\makeatother

\newcommand{\vp}{\varphi}

 \def\bE {{\mathbb E}}

 \def\bN {{\mathbb N}} 
\def\bP {{\mathbb P}}  \def\bR {{\mathbb R}}

\def\1{{\mathbbm{1}}}
\def\nn{\nonumber}

\def\qed{{\hfill $\Box$ \bigskip}}
\def\eps{\varepsilon}
\def\wt{\widetilde}
\def\vphi{\varphi}

\def\pf{\noindent{\bf Proof.} }

\begin{document}

\title[Schauder estimates]{Schauder estimates in generalized H\"older 
spaces}
\author{Jongchun Bae}
\email{bjc0204@snu.ac.kr}
\address{Department of Mathematical Sciences, Seoul National University, 
Building 27, 1 Gwanak-ro, Gwanak-gu
Seoul 151-747, Republic of Korea}

\author{Moritz Kassmann}
\email{moritz.kassmann@uni-bielefeld.de}
\address{Universit\"{a}t Bielefeld, Fakult\"at f\"ur Mathematik, Postfach 
100131, D-33501 Bielefeld, Germany}

\thanks{Support by the German Science Foundation (DFG) and the National 
Research Foundation of Korea (NRF), supporting
the initiation of an international collaboration (RO 1195/10 of DFG and 
2012K2A5A6047864 of
NRF), is gratefully acknowledged.}

% AMS subject classifications (used in AMS journals)

% AMS keywords (used in AMS journals)

\begin{abstract}
We prove Schauder estimates in generalized H\"older spaces $C^\psi(\R^d)$. 
These spaces are characterized by a general modulus of continuity $\psi$, which 
cannot be represented by a real number.  We consider linear operators $\cL$ 
between such 
spaces. The operators $\cL$ under consideration are integrodifferential 
operators with a functional  order of differentiability $\vp$ which, again, 
is not represented by a real number. Assuming that $\cL$ has 
$\psi$-continuous coefficients, we prove that solutions $u \in 
C^{\vp\psi}(\R^d)$ to linear equations $\cL u = f \in 
C^\psi(\R^d)$ satisfy a priori estimates in $C^{\vp\psi}(\R^d)$.
\end{abstract}

\maketitle

\section{Introduction}\label{sec:intro}

Schauder estimates are a central tool in the study of classical solutions to 
differential equations with H\"older continuous coefficients. In short, the 
idea of this approach is to view these equations on a small scale as 
perturbations to equations with constant coefficients. This approach 
allows to use ideas from potential theory when treating equations with variable 
coefficients. An exposition of this method can be found in any serious textbook 
on partial differential equations. 

For differential operators of second order, the main assertion in this field 
is an estimate of the form 
\[ \|u\|_{C^{2+\alpha}} \leq c \left( \|u\|_{C^{0}} + \|f\|_{C^{\alpha}} \right) \]
for all solutions $u$ to elliptic equations of second order $L u = f$ 
\cite[Theorem 6.2]{GiTr83} with some positive generic constant $c$, depending 
on the ellipticity of $L$, the dimension $d$ and the number $\alpha \in (0,1)$. 
The order $2$ and the value $\alpha$ are independent quantities in this 
estimate. The estimate holds for different values of $\alpha$ but, as shown in 
\cite{Bas09, DoKi13, RoSe14, JiXi14}, it holds analogously for solutions to 
integrodifferential equations where the driving operator is an 
integrodifferential operator with fractional order of differentiability $\beta 
\in (0,2)$. The prototype of such an operator is provided by the 
fractional Laplace operator $(-\Delta)^{\beta/2}$, which can be defined for $u 
\in 
C^2_c(\R^d)$ by 
\begin{align}\label{eq:defi_frac-laplace}
 - (-\Delta)^{\beta/2} u (x) = \cC_{\beta,d} \int\limits_{\R^d} \big( u(x+h) - 
2 u(x) + u(x-h) \big) |h|^{-d-\beta} \d h \,,
\end{align}
where $\cC_{\beta,d}$ is an appropriate constant. Note that 
$\cF((-\Delta)^{\beta/2} u)(\xi) = |\xi|^\beta \cF(u)(\xi)$, where $\cF$ 
denotes the Fourier transform. Hence the name is ``fractional Laplace operator''. 

The aim of this article is to prove Schauder estimates for a finer scale of 
function spaces and of operators at the same time. Let us explain this idea 
step by step. The H\"older space $C^\alpha(\R^d)$ is characterized by the 
number 
$\alpha \in (0,1)$ which appears in the bound of the modulus of continuity:
\begin{align*}
v \in C^\alpha(\R^d) \quad \Leftrightarrow \quad \sup_{x \in 
\R^d} |v(x)| + \sup_{x,y \in 
\R^d}\frac{|v(x)-v(y)|}{|x-y|^\alpha} < + \infty
\end{align*}
We will study estimates in more general space $C^\psi(\R^d)$ where a
function 
$\psi: (0,1] \to (0,\infty)$
 is used to replace $|x-y|^\alpha$ in 
the above expression by 
$\psi(|x-y|)$. In this sense, we study a much finer scale of function spaces. 
This scale is of particular interest when studying mapping properties 
of integrodifferential operators because for them such scales 
turn out to be natural. Note that generalized H\"older spaces have been studied 
in various settings and for 
very long. We mention several articles in \autoref{sec:spaces} when we define 
and discuss these spaces. 

Schauder estimates for integrodifferential operators 
are proved in \cite{Bas09, DoKi13, 
RoSe14, JiXi14, KiKi15} in different contexts. Our approach is inspired by the 
straightforward approach in \cite{Bas09}. We also make use of recent 
developments in potential theory obtained in \cite{KKK13} 
when studying the translation invariant case. The new contribution of this work 
is twofold. On the one hand, we allow the integration kernels to have general 
singularities at $h=0$. On the other hand, we study the resulting 
a priori estimates in a much finer scale of function spaces. Although these 
developments could be approached separately, the main new finding of our work 
is that they naturally belong together. This phenomenon does not exist in 
the study of differential equations of second order.

Let us first discuss the function spaces $C^\psi (\R^d)$. 
We assume $\lim_{r\to 0+} \psi(r) = 0$.
%, $\psi(r)\ge \psi(1)$ for all $r>1$ and that $\psi$ stays away from zero on compact subsets of $(0,\infty)$. 
In order to describe the order of 
differentiability induced by a particular function $\psi$, we need to 
introduce two indices. We define indices $M_\psi$ and $m_\psi$ by  
\begin{align}
M_\psi &= \inf \big\{ \alpha\in \bR | r \to \tfrac{\psi(r)}{r^\alpha}  \text{ 
is almost  decreasing in } (0,1] \big\} \,, \label{def:M_psi}\\
m_\psi &= \sup \big\{ \alpha\in \bR | r \to \tfrac{\psi(r)}{r^\alpha}  \text{ 
is almost 
increasing in } (0,1] \big\} \,.\label{def:m_psi}
\end{align}
See the definition for the almost monotonicity in \autoref{sec:spaces}.
Note that if $\psi$ is a regularly varying function
of order $\alpha \in (0,1)$ at zero like $\psi(r)=r^\alpha$ or 
$\psi(r)=r^\alpha |\ln(\frac2r)|$ we find $M_\psi = m_\psi = \alpha$. We denote the closed interval $[m_\psi,M_\psi]$ by $I_\psi$. This interval describes the range of 
the functional order of differentiability induced by $\psi$. For example, the 
condition $I_\psi \subset (0,1)$ implies that 
%for $\alpha = m_\psi$ und $\beta = M_\psi$ we know
 $C^{M_\psi}(\R^d) \subset C^\psi(\R^d) \subset 
C^{m_\psi}(\R^d)$. On the other hand, cases $I_\psi \cap \N \ne \emptyset$ lead 
to well-known technical difficulties which we want to avoid. See 
\autoref{sec:spaces} for a more detailed discussion of the spaces 
$C^\psi(\R^d)$ including appropriate norms.

Our ultimate goal is to study integrodifferential operators which are not 
translation invariant, i.e., which have state dependent kernels. As it is 
usually done in the theory of 
Schauder estimates, we first study the translation invariant case, i.e., we 
study integrodifferential operators with constant coefficients. Since, in this case, our 
assumptions imply that these operators satisfy the maximum principle 
and generate L\'evy processes, we can employ techniques from potential theory.

Let us define the integrodifferential operators under consideration. 
Let $\vphi:(0,\infty) \to 
(0,\infty)$ be a smooth function with $\vphi(1)=1$. We assume that 
the function $\phi$ defined by $\phi(r) = 
\vphi(r^{-1/2})^{-1}$ is a Bernstein function, i.e., satisfies $(-1)^n 
\phi^{(n)}(r) \le 0$ for every $n \in \N$. Furthermore, we assume the scaling 
condition
\begin{align}
a_1\lambda^{\delta_1} \phi(r) \le \phi(\lambda r) \le a_2 \lambda^{\delta_2} 
\phi(r) \qquad (\lambda \ge 1, r\in (0,\infty)) \label{eq:scaling-cond-phi}
\intertext{or, equivalently,}
a_1\lambda^{2\delta_1} \vphi(r) \le \vphi(\lambda r) \le a_2 \lambda^{2\delta_2} 
\vphi(r) \qquad (\lambda \ge 1, r\in (0,\infty)) \label{eq:scaling-cond-vphi}
\end{align}
for some constants $0<\delta_1\le \delta_2< 1$, $a_1\in(0,1]$, and $a_2\in 
[1,\infty)$. Typical examples are given by $\vphi(s) = s^\alpha$ or 
$\vphi(s) = s^\alpha \ln(1+s^\beta)$ for $\alpha, \beta, \alpha + \beta \in 
(0,2)$. In particular, we point out $I_\vphi\subset [2\delta_1, 2\delta_2]$. 

Let $a_0$ be a measurable function on $\bR^d \setminus \{0\}$ into $[\Lambda_1, \Lambda_2]$ for some positive 
constants $\Lambda_1$, $\Lambda_2$. 
In the case $M_\vphi \in [1,2)$ we define $\cL_0$ by
\begin{align}\label{def:L_0}
\begin{split}
\mathcal{L}_0 u(x) &= \int_{\bR^d} (u(x+h)-u(x)-\nabla u(x)\cdot h 
\1_{\{|h|\le 1\}} 
) \frac{a_0(h)}{|h|^{d}\vphi(|h|)} \d h,
\end{split}
\end{align}
and in the case $M_\vphi \in (0,1)$ by
\begin{align}\label{def:L_0-no-grad}
\begin{split}
\mathcal{L}_0 u(x) &= \int_{\bR^d} (u(x+h)-u(x)) 
\frac{a_0(h)}{|h|^{d}\vphi(|h|)} \d h
\end{split}
\end{align}
for all continuous functions $u:\R^d \to \R$ for which the integral is 
well defined for every point $x \in \R^d$. Note that this domain 
$\cD(\cL_0)$ includes functions $u\in 
C^2(\R^d)$ which we assume to be bounded. The class of operators 
$\mathcal{L}_0$ 
is significantly larger than those of \eqref{eq:defi_frac-laplace}. The main 
difference is that the order of differentiability of $\mathcal{L}_0$ is 
represented by the 
function $\vphi$ and cannot be represented by a single number. Note that 
$\vphi$ may be chosen as not regularly varying at zero.  As we will see, our 
scale of function spaces $C^\psi$ is well suited to formulate mapping 
properties of such operators.

Our Schauder estimate for translation invariant operators reads as follows. 
\begin{theorem}\label{theo:mainL_0-intro} Let $\vphi$ and $\psi$ be functions described above. Suppose $I_\psi\subset(0,1)$ and 
$I_{\vphi\psi}\subset(0,1)\cup (1,2)\cup 
(2,3)$. There 
exists $C_1$ such that for 
$u \in C^{\vphi \psi}(\R^d)\cap C_\infty(\R^d)$ and  $f \in C^\psi(\R^d)$ 
satisfying $\cL_0 u=f$, 
the following estimate holds:
\begin{align*}
\| u\|_{C^{\vphi\psi}} \le C_1( \|u\|_{C^0} + \| f\|_{C^{\psi}}).
\end{align*}
\end{theorem}
We prove this result in \autoref{sec:translation-invariant-operators} using a 
semigroup approach. In our proofs we benefit from ideas in \cite{Bas09} and 
\cite{KKK13}. Once \autoref{theo:mainL_0-intro} is established, we can use a 
perturbation 
argument to treat integrodifferential operators with variable coefficients $a: 
\R^d\times \R^d \to [\Lambda_1, \Lambda_2]$. Let $\cL$ be defined by 
\begin{align}\label{def:L1+}
\cL u(x) = \int_{\bR^d} \left(u(x+h)-u(x)-\nabla u(x) \cdot h 
\1_{\{|h|\le 
1\}}\right) \frac{a(x,h)}{|h|^d\vphi(|h|)} \d h
\end{align}
when $M_\vphi \in[1,2)$ and by \begin{align}\label{def:L1-}
\cL u(x) = \int_{\bR^d} 
\left(u(x+h)-u(x)\right)\frac{a(x,h)}{|h|^d\vphi(|h|)} \d h
\end{align}
when $M_\vphi\in(0,1)$ for all continuous functions $u:\R^d \to \R$ for which 
the integral is well defined for every point $x \in \R^d$. Note that 
this domain equals $\cD(\cL_0)$. The coefficient function 
$a:\R^d \times \R^d \to [\Lambda_1, \Lambda_2]$ is 
assumed to satisfy 
\begin{alignat}{2}
\sup_{x\in\bR^d}\sup_{|h|>0}|a(x+z,h)-a(x,h)| &\le \Lambda_3 \psi(|z|)
\qquad && (|z|\le 1) \label{eq:assum_axh-cont}
\end{alignat}
for some positive constant $\Lambda_3 \geq 1$. 
%Here, the function $\psi$ is 
%assumed to satisfy the assumptions above. 
This condition requires the function $x 
\mapsto a(x,h)$ to be $\psi$-continuous uniformly in $h$.

We have already mentioned that the definition of H\"older and H\"older-Zygmund 
spaces is delicate when the order of differentiability is an 
integer. In order to formulate our main result we need to exclude this case. We 
assume further:
\begin{align}\label{ass:I_contain_no_integer}
[m_{\vphi\psi}, M_{\vphi\psi}] \cap \N = \emptyset\,, \qquad 
M_\vphi \vee M_\psi < m_{\vphi \psi} \,.
\end{align}
% We define a partial order in the class of intervals $I_\psi$'s as  $I_{\psi_1} 
% <I_{\psi_2}$ if  $M_{\psi_1} <m_{\psi_2}$.
Let us formulate the main result of this work.

\begin{theorem}\label{theo:main}
Assume that, in addition to the assumptions of 
\autoref{theo:mainL_0-intro}, condition   
\eqref{ass:I_contain_no_integer} is satisfied.  In the case $1 \in I_\vphi$, we
assume $a(x,h)=a(x,-h)$ for all $x,h\in\R^d$. 
Then there exists a positive constant $C_2$ such that for every $f \in C^{ 
\psi}(\R^d)$ and every solution $u \in C^{\vphi\psi}(\bR^d)$ to 
the equation $\cL u =f$ the following estimate holds:
\begin{align*}
\| u\|_{C^{\vphi\psi}} \le C_2 \big( \|u\|_{C^0} +  \| f\|_{C^{ \psi}} \big) 
\,.
\end{align*}
\end{theorem}

Let us make a few comments. Note that \autoref{theo:main} trivially implies 
\autoref{theo:mainL_0-intro}. The assumption  
$a(x,h)=a(x,-h)$ for all $x,h\in\R^d$ in the case $1 \in I_\vphi$ is natural 
due to the appearance of the gradient term in the definition of $\cL$. Note 
that this assumption needs to be added to \cite[Corollary 
5.2]{Bas09} in order for the corollary to be correct. The first part 
of \eqref{ass:I_contain_no_integer} is 
natural and resembles the fact that Lipschitz function space is not equal to the space $C^1(\R^d)$. The other parts of 
\eqref{ass:I_contain_no_integer}  would 
vanish if we restricted ourselves to the (large) class of regularly varying 
functions $\psi$ and $\vphi$.  

It is important to note that the a priori estimate provided by 
\autoref{theo:main} is the best possible. This follows from the mapping 
properties of $\cL$, which are given in the following theorem. We defer the 
proof of this result to \autoref{sec:maprop}.
\begin{theorem}\label{theo:maprop}
Suppose that $M_\vphi \vee M_\psi < m_\vphi + m_\psi$. Assume $\cL$ and 
$a(\cdot,\cdot)$ satisfy \eqref{def:L1+} resp. \eqref{def:L1-}, and 
\eqref{eq:assum_axh-cont}. Furthermore, if $1\in I_\vphi$ we assume that 
$a(x,h)= a(x,-h)$ for every $x,h \in \R^d$. Then the operator $\cL$ is a 
continuous operator from $C^{\vphi\psi}(\R^d)$ to $C^{\psi}(\R^d)$.
\end{theorem}

The article is organized as follows. In \autoref{sec:spaces} we define and 
study the generalized H\"older spaces $C^\psi(\R^d)$. The proof of 
\autoref{theo:main} relies on a perturbation technique. First, we study the 
operator $\cL_0$ which is obtained by ``freezing'' the coefficients via $a_0(h) 
= a(x_0,h)$ for some arbitrary but fixed point $x_0 \in \R^d$. In 
\autoref{sec:translation-invariant-operators} we derive estimates on the 
transition density for the semigroup generated by $\cL_0$. 
\autoref{sec:translation-invariant-operators} also contains the proof of 
\autoref{theo:mainL_0-intro}. In \autoref{sec:main-proof} we prove our 
main result, \autoref{theo:main}. The proof of \autoref{theo:maprop} is 
given in \autoref{sec:maprop}.

%%%%%%%%%%%%%%%%%%%%%%%%%%%%%%%%%%%%%%%%%%%%%%%%%%%%%%%%%%%%%%%%%%%%%%%%%%%%
%%%%%%%%%%%%%%%%%%%%%%%%%%%%%%%%%%%%%%%%%%%%%%%%%%%%%%%%%%%%%%%%%%%%%%%%%%%%
\section{Generalized H\"older spaces}\label{sec:spaces}
%%%%%%%%%%%%%%%%%%%%%%%%%%%%%%%%%%%%%%%%%%%%%%%%%%%%%%%%%%%%%%%%%%%%%%%%%%%%
%%%%%%%%%%%%%%%%%%%%%%%%%%%%%%%%%%%%%%%%%%%%%%%%%%%%%%%%%%%%%%%%%%%%%%%%%%%%

In this section we define the function spaces $C^\psi(\R^d)$ and discuss several 
of their properties. Unfortunately, we are not able to use results from the 
literature despite an intensive search. Since generalized smoothness of 
functions and related function spaces have been studied for several decades, it 
is likely that the results of this section have been proved 
somewhere else. Let us mention only a very few expositions which might be 
valuable for the interested reader. A very early source is \cite{BaSt56}. 
Several cases and results are established in \cite{Ulj68, Gol72, KiLi87}. Some 
more recent works include \cite{FaLe06, KrNi12} where many more references 
can be 
found.      

We denote by  $C^0(\R^d)$ the Banach space of real-valued, 
bounded, and continuous 
functions on $\R^d$ equipped with the norm $\| f\|_{C^0(\R^d)} = \| 
f\|_{C^0} = \sup_{x\in \R^d} |f(x)| <\infty $. For $m \in \N$ we denote by $C^m(\R^d)$ the Banach space 
of functions $f \in C^0(\R^d)$ with all derivatives $D^\gamma f \in 
C^0(\R^d)$ for $|\gamma| \leq m$. Here, we denote by $D^\gamma f$ the partial derivative $\partial_{x_1}^{\gamma_1}\cdots \partial_{x_d}^{\gamma_d} f$ and $|\gamma| = \sum_{i=1}^d \gamma_i$ for the multi index $\gamma =(\gamma_1,\ldots, \gamma_d) \in \N^d$. By 
$C(\R^d)$ we denote the Fr\'{e}chet space of real-valued 
continuous functions on $\R^d$.

Let $\psi$ be a positive real valued function on $(0,1]$ with $\psi(1)=1$ and $\lim_{r\to 0+} \psi(r) =0$.
For $j \in \N_0$ we define a seminorm $[f]_{C^{-j;\psi}(\R^d)}$ 
by
\begin{align*}
\left[f\right]_{C^{-j;\psi}(\R^d)} =\sup_{x\in\R^d}\sup_{0<|h|\le 1} 
\frac{|f(x+h)-f(x)|}{\psi(|h|)|h|^{-j}},
\end{align*}
and a vector space of functions $C^{-j;\psi}(\R^d)$ by
\begin{align*}
C^{-j;\psi}(\R^d) = \left\{ f\in C(\R^d) \big|  [f]_{C^{-j;\psi}(\R^d)}<\infty 
\right\}\,.
\end{align*}

We abuse the notation $[u]_{C^{0;\psi}} = [u]_{C^\psi}$ for the convenience.
Following 
\cite{BGT89}, for a subinterval $I$ on $(0,\infty)$ we call a function  $\psi: 
I 
\to (0,\infty)$ \emph{almost increasing} 
if there is a constant $c\in (0,1]$ such that 
$c \psi(r) \le \psi(R)$ for $ r, R \in I, r \le R$. On the other hand, we call 
such 
$\psi$ \emph{almost decreasing} if there is $C\in [1,\infty)$ such that $ 
\psi(R)\le C 
\psi(r) $ for $r,R \in I$ and $ r\le R$. Recall the definition of the indices $M_\psi$ 
and $m_\psi$ from \eqref{def:M_psi} and 
\eqref{def:m_psi}. Now we can finally define the function spaces 
$C^{\psi}(\R^d)$. 

\begin{definition}\label{def:psi-spaces}
In the case of 
$m_\psi \in(0,1]$, let $C^{\psi}(\R^d)$ be defined by
\begin{align*}
C^\psi(\R^d) = \{ f\in C(\R^d) | f \in C^0(\R^d) \text{ and } 
[f]_{C^{0;\psi}}<\infty\}. 
\end{align*}
In the case of $m_\psi \in (k,k+1]$ for some $k\in \bN$, let $C^{\psi}(\R^d)$ be 
defined by 
%\begin{align}
%C^\psi(\R^d) = \{ f\in C(\R^d) | f, Df, \ldots, D^kf \in C^0(\R^d) \text{ and } D^kf \in 
%C^{-k;\psi}(\R^d) \}\,.
%\end{align}
\begin{align*}
C^\psi(\R^d) = \{ f\in C(\R^d)  \,|\,  & D^\gamma f \in C^0(\R^d) \text{ for all } 0\le |\gamma|\le k\nn \\
& \text{ and } D^\gamma f \in 
C^{-k;\psi}(\R^d) \text{ for all }  |\gamma| = k\}\,.
\end{align*}

In the case of $m_\psi \in (k,k+1]$ for some $k\in \N_0$, the $\psi$-H\"older norm $\| 
\cdot \|_{C^\psi}$ is defined by 
\begin{align*} 
\left\| f\right\|_{C^{\psi}} = \sum_{j=0}^k \left\| D^j f\right \|_{C^0} 
+ [D^kf ]_{C^{-k;\psi}(\R^d)} \,.
\end{align*}
\end{definition}
Here, we use the notation $D^0f = f$ and denote the maximum of $C^0$-norms and $C^{-k;\psi}$-seminorms of all $k$-th derivatives of $f$ by $\| D^k f\|_{C^0}$ and $[D^kf]_{C^{-k;\psi}}$ respectively. In particular, when $k=1$, we omit the exponent for the sake of brevity.

If there is no ambiguity, then we write $C^{\psi}$ instead of $C^{\psi}(\R^d)$. Let us 
start with some observations. Note that for $\alpha \in (0,1)$ the two 
seminorms 
\begin{align*}
[f]^{(1)}_{C^\alpha} = \sup_{x\in\R^d}\sup_{|h|>0} 
\frac{|f(x+h)-f(x)|}{|h|^\alpha} \quad \text{ and } \quad [f]^{(2)}_{C^\alpha} 
= \sup_{x\in\R^d}\sup_{|h|>0} \frac{|f(x+h)-2f(x)+f(x-h)|}{|h|^\alpha} 
\end{align*}
are equivalent. We prove an analogous property in our more general function 
spaces. The condition $\alpha \in (0,1)$ translates to $I_\psi \subset(0,1)$ in our 
setting. To shorten notation, let us write first-order and second-order differences as 
follows: 
\begin{align*}
\Delta_hf(x) = f(x+h)-f(x) \quad \text{ and } \quad \Delta^2_hf(x) = f(x+h)-2f(x)+f(x-h) 
\,.
\end{align*}

For the sake of brevity we use the notation 
$$ [[f]]_{C^{\psi}} =  \sup_{x\in\R^d} \sup_{0<|h|\le 1} \frac{|f(x+h)-2f(x)+f(x-h)|}{\psi(|h|)}.$$
Triangle inequality gives the trivial inequality
$$ [[f]]_{C^\psi} \le 2[f]_{C^\psi}.$$
We will show in the following lemma that the seminorm $[f]_{C^\psi}$ is bounded above by the sum of $\|f\|_{C^0}$ and a seminorm $[[f]]_{C^\psi}$. 
Summing up we get the equivalence between the two norms $\|f\|_{C^\psi}$ and $\|f\|_{C^0} + [[f]]_{C^\psi}$.
\begin{lemma}\label{lem:psi_hoelder1}
Let $I_\psi \subset(0,1)$ and $f\in C^{\psi}$. There exists a constant $C= C(\psi)$ such that
\begin{align*}
[f]_{C^\psi} \le C\left(\|f\|_{C^0} + [[f]]_{C^\psi}\right).
\end{align*}
\end{lemma}
\pf From the definition of $m_\psi$ and $M_\psi$ we choose constants $c_1\in (0,1]$ and $c_2\in[1,\infty)$  such that 
\begin{align}\label{e:pf_psi_ratio}
c_1 \left(\frac{R}{r}\right)^{m_\psi/2} \le \frac{\psi(R)}{\psi(r)} \le c_2 \left(\frac{R}{r}\right)^{(M_\psi+1)/2} \quad \text{ for } r\le R\le 1.
\end{align}
 Let $n$ be an integer greater than $(2c_2)^{2/(1-M_\psi)}$. For every $0<|h|\le 1$, we have
\begin{align*}
n\Delta_{h/n}f(x) = \Delta_hf(x) - \sum_{k=1}^{n-1} (n-k)\Delta_{h/n}^2 f(x+\frac{k}{n}h).
\end{align*}
Dividing by $\psi(|h|)$, we obtain 
\begin{align*}
\frac{n\psi(|h|/n)}{\psi(|h|)} \frac{|\Delta_{h/n}f(x)|}{\psi(|h|/n)} \le \frac{|\Delta_h f(x)|}{\psi(|h|)} +\sum_{k=1}^{n-1} (n-k)\frac{|\Delta_{h/n}^2f(x+\frac{k}{n}h)|}{\psi(|h|/n)} \frac{\psi(|h|/n)}{\psi(|h|)}.
\end{align*} 
Using \eqref{e:pf_psi_ratio} with $R=|h|$ and $r=|h|/n$ and taking supremum over $x\in \R^d$ and $0<|h|\le 1$, we get
\begin{align*}
c_2^{-1}n^{(1-M_\psi)/2} \sup_{x\in\R^d}\sup_{0<|h|<1/n} \frac{|\Delta_{h}f(x)|}{\psi(|h|)} \le [f]_{C^\psi} + \frac{(n-1)n^{1-m_\psi/2}}{2c_1}[[f]]_{C^\psi}.
\end{align*}
Since $\psi(|h|) \ge c_2^{-1}|h|^{(M_\psi+1)/2} \ge c_2^{-1}n^{-(M_\psi+1)/2}$ for $1/n \le |h|\le 1$, we also have
\begin{align*}
c_2^{-1}n^{(1-M_\psi)/2} \sup_{x\in\R^d}\sup_{1/n\le |h|\le 1} \frac{|\Delta_{h}f(x)|}{\psi(|h|)} \le 2n\|f\|_{C^0}. 
\end{align*}
Therefore, from the choice of $n$ we obtain
\begin{align*}
[f]_{C^\psi} \le  2n\|f\|_{C^0}+ (2c_1)^{-1}(n-1)n^{1-m_\psi/2} [[f]]_{C^\psi}, 
\end{align*}
which implies the result. \qed

This equivalence is allowed for the case $I_\psi \subset (1,2)$ by the following lemmas.

\begin{lemma}\label{lem:interpolation1}
If $I_{\psi} \subset (1,2)$, then for small $\eps >0$ there exists a constant 
$C=C(d,\psi,\eps)>0$  such that
\begin{align}\label{e:Df_eps1}
\| D f\|_{C^0} \le C\|f \|_{C^0} + \eps [ D f ]_{C^{-1;\psi}}. 
\end{align}
If  $I_{\psi} \subset (2,3)$, then for small $\eps >0$ there exists a constant 
$C=C(d,\psi,\eps)>0$ such that
\begin{align}\label{e:Df_eps2}
\|D f\|_{C^0}+ \| D^2 f\|_{C^0} \le C\|f \|_{C^0} + \eps [ D^2 f ]_{C^{-2;\psi}}. 
\end{align}
\end{lemma}
\pf First consider the case $I_\psi \subset(1,2)$. 
Fix $1\le i \le d$ and $x$ be a point in $\bR^d$. Let $c_1\in(0,1]$ such that 
\begin{align*}
c_1\left(\frac{R}{r}\right)^{(m_\psi+1)/2} \le \frac{\psi(R)}{\psi(r)},  \quad \text{ for } r \le R \le 1.
\end{align*}
The case when $[ D_if ]_{C^{-1;\psi}} = 0$ is trivial, so we suppose 
not. Define 
$N=(\|f\|_{C^0} / [D_if]_{C^{-1;\psi}})^{2/(m_{\psi}+1)}$ if $\|f\|_{C^0}\le [D_if]_{C^{-1;\psi}}$ and $N=(\|f\|_{C^0} / [D_if]_{C^{-1;\psi}})^{2/(M_{\psi}+2)}$ otherwise. We may only consider the first case because the proof for the other case is the same. By the mean value 
theorem, 
there exists $x'$ on the line segment between $x$ and $x + Ne_i$ such that
$$ | D_i f(x') | = \frac{|f(x+N e_i)-f(x)|}{N} \le 
\frac{2\|f\|_{C^0}}{N}.$$
Thus
$$|D_i f(x)| \le |D_if(x')| + |D_if(x')-D_if(x)| \le \frac{2\|f\|_{C^0}}{N} 
+ 
c_1^{-1}[D_if]_{C^{-1;\psi}} \psi(N)N^{-1}.$$ 
With the fact $\psi(N) \le c_1^{-1} N^{(m_{\psi}+1)/2}$ and the choice of $N$,
\begin{align}\label{e:aux1_Df}
 |D_i f(x)| \le (2+c_1^{-2})\|f\|_{C^0}^{1-2/(m_{\psi}+1)} 
[D_if]_{C^{-1;\psi}}^{2/(m_\psi+1)}.
\end{align}
Taking the supremum over $x\in \bR^d$ and then applying the inequality
\begin{align*}
 r^{\theta} s^{1-\theta} \le r + s, \quad r,s >0, \theta \in (0,1),
\end{align*}
we obtain
$$ \| D_if \|_{C^0} \le (2+c_1^{-2}) ( \| f\|_{C^0} +  [ D_if ]_{C^{-1;\psi}}).$$ By the scaling argument we get \eqref{e:Df_eps1}.

Now, we assume $I_\psi\subset(2,3)$. Let $c_2$ be a constant such that 
\begin{align*}
c_2\left(\frac{R}{r}\right)^{(m_\psi +2)/2} \le \frac{\psi(R)}{\psi(r)}, \quad \text{ for } r \le R \le 1.
\end{align*}
 By the above argument used to obtain \eqref{e:aux1_Df}, we have
\begin{align*}
\|D_if\|_{C^0} \le 3\|f\|_{C^0}^{1/2}\|D_{ii}f\|_{C^0}^{1/2} 
\end{align*}
Thus it suffices to show the second result for the left hand side replaced by $\|D^2 f\|_{C^0}$. The same argument for \eqref{e:aux1_Df} we get
\begin{align*}
|D_{ij} f(x) | &\le (2+c_2^{-2})\|D_j f\|_{C^0}^{1-2/m_\psi} [D_{ij} f]_{C^{-2;\psi}}^{2/m_\psi}\\
&\le 3(2+c_2^{-2})\|f\|_{C^0}^{1/2-1/m_\psi} \|D^2f\|_{C^0}^{1/2-1/m_\psi}[D_{ij} f]_{C^{-2;\psi}}^{2/m_\psi},
\end{align*}
which implies
\begin{align*}
\|D^2 f\|_{C^0} \le c_3\|f\|_{C^0}^{(m_\psi-2)/(m_\psi+2)} [D^2f ]_{C^{-2;\psi}}^{4/(m_\psi+2)} \le c_3 (\|f\|_{C^0} + [D^2 f]_{C^{-2;\psi}}).
\end{align*}
By the scaling argument we get \eqref{e:Df_eps2}. \qed

\begin{lemma}\label{lem:psi_hoelder2}
Assume $I_\psi\subset(1,2)$ and $f\in C^{\psi}(\R^d)$.  Then there exists a constant $C=C(\psi)$ such that
\begin{align}\label{e:proof-hoelder2}
[D_i f]_{C^{-1;\psi}} \le C\left(\|f\|_{C^0} + [[f]]_{C^\psi}\right),
\end{align}
for every $i=1,\ldots, d$.
\end{lemma}
{\bf Proof:} Define $\bar{\psi}(r): 
= r^{-1}\psi(r)$. Then $I_{\bar{\psi}} \subset (0,1)$. First note that  it is shown in \autoref{lem:psi_hoelder1}  that the seminorm  $[\ \cdot\ ]_{C^{\bar{\psi}}}(=[\ \cdot\ ]_{C^{-1;\psi}})$ is bounded by $\|\cdot \|_{C^0} + [[\ \cdot\ ]]_{C^{\bar{\psi}}}$. Choose  $c_1 \geq 1$ such that for every $g\in C^{-1;\psi}$
\begin{align}\label{e:proof-psi-onediff<twodiff}
[g]_{C^{-1;\psi}} \le c_1 \left(\|g\|_{C^0} + [[g]]_{C^{\bar{\psi}}}\right)\,.
\end{align}
Since $D_if\in C^{-1;\psi} $ for any $f\in C^\psi$, 
\begin{align*}
[D_if]_{C^{-1;\psi}} \le c_1\left(\|D_if\|_{C^0} + [[D_if]]_{C^{\bar{\psi}}}\right).
\end{align*}
We claim that 
\begin{align}\label{e:proof-psi-hoelder-Df<f}
[[D_if]]_{C^{\bar{\psi}}} \le c_2 (\|D_i f\|_{C^0} + [[f]]_{C^\psi}), \quad i=1,\ldots, d
\end{align}
for some constant $c_2$ not depending on $f$.
If we prove the claim, then the following estimate from \autoref{lem:interpolation1} below
$$\|D_if\|_{C^0} \le  (2c_1(1+c_2))^{-1} [D_if]_{C^{-1;\psi}}+ c_3\|f\|_{C^0} $$
implies  \eqref{e:proof-hoelder2} with $C= 2c_1(1+c_2)(1+c_3)$ .

In order to prove \eqref{e:proof-psi-hoelder-Df<f} we only 
consider the case $i=1$. The remaining cases can be dealt with analogously. For $k,h 
\in \bR^d$  we obtain
\begin{align*}
&\quad |k||\Delta_h^2(D_1f)(x)|\\
&=\bigg| |k|D_1f(x+h)-\Delta_{|k|e_1}f(x+h) + \Delta_{|k|e_1}f(x+h)\\
&\quad -2\big( |k|D_1f(x)-\Delta_{|k|e_1}f(x) \big) - 2\Delta_{|k|e_1}f(x)\\
&\quad +|k|D_1f(x-h)-\Delta_{|k|e_1}f(x-h) + \Delta_{|k|e_1}f(x-h)\bigg|\\
&\le |k|\big|D_1f(x+h)-D_1f(x+h+\theta_1|k|e_1)\big| + 
2|k| \big| D_1f(x)-D_1f(x+\theta_2|k|e_1) \big| \\
&\quad + |k| \big|D_1f(x-h)-D_1f(x-h+\theta_3|k|e_1) \big| + |\Delta_h^2f(x+|k|e_1)| + 
|\Delta_h^2f(x)|
\end{align*}
for some real numbers $\theta_1,\theta_2,\theta_3\in[0,1]$. Let $c_4\in(0,1] $ be a constant such that
\begin{align}\label{e:pf_psi_ratio_1}
 c_4\left(\frac{R}{r}\right)^{(1+m_\psi)/2} \le \frac{\psi(R)}{\psi(r)} \quad \text{ for } r \le R \le 1.
 \end{align}
If $|k|\le |h| \le 1$ then the sum of the first three terms is bounded by
\begin{align}\label{e:proof-psi-onediff}
4c_4^{-1}[D_1f]_{C^{-1;\psi}}\psi(|k|),
\end{align}
and the sum of the last two terms is bounded by 
\begin{align}\label{e:proof-psi-twodiff}
2[[f]]_{C^{\psi}}\psi(|h|).
\end{align}
Using \eqref{e:proof-psi-onediff<twodiff}, \eqref{e:proof-psi-onediff} and \eqref{e:proof-psi-twodiff} we have
\begin{align*}
|k||\Delta_h^2(D_1f)(x)| \le 4c_1c_4^{-1}([[D_1f]]_{C^{\bar{\psi}}} + \|D_1f\|_{C^0})\psi(|k|) + 
2[[f]]_{C^{\psi}}\psi(|h|)
\end{align*}
for $|k|\le |h|\le 1$.  
Taking $k=\eps h$ with $\eps = (8c_1c_4^{-2})^{2/(1-m_\psi)}$ gives us
\begin{align*}
& \eps |h||\Delta_h^2(D_1f)(x)| \le 4c_1c_4^{-1}([[D_1f]]_{C^{\bar{\psi}}} + \|D_1f\|_{C^0})\psi(\eps |h|)+ 
2[[f]]_{C^{\psi}}\psi(|h|).
\end{align*}
Dividing both sides by $\eps\psi(|h|)$ and using \eqref{e:pf_psi_ratio_1}, we obtain
\begin{align}\label{e:pf-D^2_h1}
&\quad \frac{|\Delta_h^2(D_1f)(x)|}{\psi(|h|)|h|^{-1}} \le 4c_1c_4^{-2} 
\eps^{(m_\psi -1)/2}([[D_1f]]_{C^{\bar{\psi}}}+\|D_1f\|_{C^0}) + 2\eps^{-1}[[f]]_{C^\psi}
\end{align}
for $|h|\le 1$. 
Taking supremum to \eqref{e:pf-D^2_h1} over $x\in\R^d$ and $0<|h|\le 1$ we have an inequality
\begin{align*}
[[D_1f]]_{C^{\bar{\psi}}} \le 4c_1c_4^{-2}\eps^{(m_\psi-1)/2}\left([[D_1f]]_{C^{\bar{\psi}}} + \|D_1f\|_{C^0}\right) + 
2\eps^{-1}[[f]]_{C^{\psi}}
\end{align*}
which implies \eqref{e:proof-psi-hoelder-Df<f} with $c_2= 1+ 4(8c_1c_4^{-2})^{2/(m_\psi-1)}$.\qed 

By summing up the results in \autoref{lem:psi_hoelder1}, \autoref{lem:interpolation1} and \autoref{lem:psi_hoelder2} we get the following equivalence.
\begin{proposition}\label{p:equi_norms1}
Let $I_\psi \subset (0,1) \cup (1,2)$. For $f\in 
C^{\psi}$ the norm $\| f\|_{C^\psi}$ is equivalent to the norm
$$ \|f\|_{C^0} + [[f]]_{C^{\psi}}.$$
\end{proposition}

\begin{proposition}\label{p:ineq_norms}
Assume $I_{\psi_1}, I_{\psi_2}\subset(0,1)\cup (1,2)\cup (2,3)$ and $M_{\psi_1} < m_{\psi_2}$. For any $0<\eps <1$,  there exists 
$C=C(d,\psi_1,\psi_2,\eps)>0$ such that
\begin{align*}
\left\| f\right\|_{C^{\psi_1}} \le C \left\|f\right\|_{C^0} + \eps \left\| 
f\right\|_{C^{\psi_2}}
\end{align*}
\end{proposition}
\pf We first consider the case $I_{\psi_2} \subset (0,1)$. Let $c_1$ and $c_2$ be the constants such that $\psi_1(|h|) \ge c_1|h|^{(2M_{\psi_1}+m_{\psi_2})/3}$ and $\psi_2(|h|) \le 
c_2|h|^{(M_{\psi_1}+2m_{\psi_2})/3}$.
Let $h_0 = 
(c_1c_2^{-1}\eps)^{3/(m_{\psi_2}-M_{\psi_1})}$. If $|h| \le h_0$, then
$$ \frac{|\Delta_hf(x)|}{\psi_1(|h|)} \le  
[f]_{C^{0;\psi_2}} \frac{\psi_2(|h|)}{\psi_1(|h|)}\le \eps [ f ]_{C^{0;\psi_2}}  .$$
If $h_0< |h| \le 
1$, then
$$ \frac{|\Delta_hf(x)|}{\psi_1(|h|)} \le \frac{2}{\psi_1(|h|)} \|f\|_{C^0}  \le 
c_3\psi_1(h_0)^{-1} \|f\|_{C^0}.$$
Combining the above two inequalities and taking supremum, we have
\begin{align}\label{e:pf-interpolation2}
 \sup_{x\in\bR^d} \sup_{0<|h|\le 1} \frac{|\Delta_hf(x)|}{\psi_1(|h|)} \le 
c_4\|f\|_{C^0} + \eps [f]_{C^{0;\psi_2}}.
\end{align}

Now we consider the case $I_{\psi_2} \subset(1,2)$. When $I_{\psi_1} \subset (0,1)$, it follows from \eqref{e:Df_eps1} that 
$$\|f\|_{C^{\psi_1}} \le c_5(\|f\|_{C^0} + \|Df\|_{C^0}) \le c_6\|f\|_{C^0} + \eps [Df]_{C^{-1;\psi_2}}.$$ 
When $I_{\psi_1} \subset(1,2)$, it also follows from \eqref{e:Df_eps1} that 
$$ \| f\|_{C^0} + \| D f \|_{C^0} \le c_8 \| f\|_{C^0} + \eps [ D f]_{C^{-1;\psi_2}}.$$
Since $[Df]_{C^{-1;\psi}} = [Df]_{C^{0;\bar{\psi}}}$, where $\bar{\psi}(r) = r^{-1}\psi(r)$, applying \eqref{e:pf-interpolation2} to $Df$ with $\bar{\psi}_1$ and $\bar{\psi}_2$ follows that
$$[Df]_{C^{-1;\psi_1}} \le c_9 \|f\|_{C^0} + \eps [Df]_{C^{-1;\psi_2}}.$$
The remaining case $I_{\psi_2} \subset (2,3)$ is also proved by the same argument combined with  \eqref{e:Df_eps2}. \qed

The product rule of derivatives gives us the following lemma.
\begin{lemma}\label{l:norm_fg}
Assume $I_{\psi}\subset (k,k+1)$ for $k\in \bN$. Then, there exists a constant $c_1=c_1(d,k,\psi) >0$ such that
$$ \|fg\|_{C^\psi} \le c_1 \|f\|_{C^\psi}\|g\|_{C^\psi}.$$
\end{lemma}
\pf By the product rule and the fact that $[D^j f]_{C^{-k;\psi}} \le c_2\|D^{j+1}f\|_{C^0}$ for $j\le k$ we can obviously obtain the result. \qed

%%%%%%%%%%%%%%%%%%%%%%%%%%%%%%%%%%%%%%%%%%%%%%%%%%%%%%%%%%%%%%%%%%%%%%%%%%%%
%%%%%%%%%%%%%%%%%%%%%%%%%%%%%%%%%%%%%%%%%%%%%%%%%%%%%%%%%%%%%%%%%%%%%%%%%%%%
\section{The translation 
invariant case}\label{sec:translation-invariant-operators}
%%%%%%%%%%%%%%%%%%%%%%%%%%%%%%%%%%%%%%%%%%%%%%%%%%%%%%%%%%%%%%%%%%%%%%%%%%%%
%%%%%%%%%%%%%%%%%%%%%%%%%%%%%%%%%%%%%%%%%%%%%%%%%%%%%%%%%%%%%%%%%%%%%%%%%%%%

The aim of this section is to prove \autoref{theo:mainL_0-intro}. Let us recall
the main assumptions from the introduction. As explained in \eqref{def:L_0}, 
\eqref{def:L_0-no-grad} we study operators of the form
\begin{align}\label{def:operator_const_coeff}
\mathcal{L}_0 u(x) &= \int_{\bR^d} (u(x+h)-u(x)-\nabla u(x)\cdot h 
\1_{\{|h|\le 1\}} 
) \frac{a_0(h)}{|h|^{d}\vphi(|h|)} \d h,
\end{align}
where $a_0 : \bR^d \setminus \{0\} \to 
[\Lambda_1, \Lambda_2]$ is a measurable function and $\Lambda_1$, $\Lambda_2$ 
are positive numbers. The domain of this operator $\cL_0$ contains bounded 
and smooth functions, e.g., $u \in C^2(\R^d)$. Recall that we assume that 
$\vphi:(0,\infty) \to 
(0,\infty)$ is a smooth function with $\vphi(1)=1$ and  
the function $\phi$ defined by $\phi(r) = 
\vphi(r^{-1/2})^{-1}$ is a Bernstein function, i.e., satisfies $(-1)^n 
\phi^{(n)}(r) \le 0$ for every $n \in \N$. Furthermore, we assume the scaling 
condition \eqref{eq:scaling-cond-phi} or, equivalently, 
\eqref{eq:scaling-cond-vphi}.

Our main idea is to apply methods from potential theory. Note that 
\[ \nu(\!\d h) = \frac{a_0(h)}{|h|^d\vphi(|h|)} \d h \]
defines a L\'evy measure with respect to a centering function $\1_{\{|h|\le 1\}}$. This measure $\nu$ induces a strongly 
continuous contraction semigroup $(P_t)$ on the Banach space $C_\infty(\R^d)$ 
of continuous functions from $\R^d$ to $\R$ that vanish at infinity. We 
write $C_\infty$ instead of $C_\infty(\R^d)$. In 
fact, $(P_t)$ is also a semigroup on $C^0$ but not strongly continuous in 
general. Denote by $C^2_\infty=C^2_\infty(\R^d)$ the space of 
functions from 
$C_\infty$ with the property that all derivatives up to order $2$ 
are elements from $C_\infty$. The infinitesimal generator $(A, \cD(A))$ of 
the semigroup $(P_t)$ satisfies $A u = \cL_0 u$ for every $u \in 
C^2_\infty$. 
%and $\cL_0$ as in \eqref{def:operator_const_coeff} when $M_\vphi 
%\ge 1$, and without the $\nabla u(x)$ term when $M_\vphi <1$. 

Our aim is to study the semigroup $(P_t)$. To do this, we first consider a
subordinate Brownian motion $X$ with subordinator whose Laplace exponent is 
$\phi$, see \autoref{subsec:subord-BM}. If we denote by $(Q_t)$ the semigroup 
of $X$, i.e., if 
$$ Q_t f(x) = \int_{\bR^d} q_d(t,x-y) f(y) dy,$$
for $f \in C_\infty(\R^d)$, then its infinitesimal 
generator $(L, \cD(L))$ acts on functions $f \in C^2_\infty(\R^d)$ in the 
following form:
\begin{align}\label{eq:generator-sub-BM}
Lf(x) = \int_{\R^d} (f(x+h)-f(x)-\nabla f(x)\cdot h \1_{\{|h|\le 1\}})J(h) \d h.
\end{align}

The values of the so-called jumping function $J(h)$ are known to be 
comparable  to $|h|^{-d}\phi(|h|^{-2})$. Since $\vphi(r)=\phi(r^{-2})^{-1}$, 
these values are also comparable to $\frac{a_0(h)}{|h|^{d}\vphi(|h|)}$ 
appearing in the definition of $\cL_0$ in \eqref{def:operator_const_coeff}. 
That is why estimates of the semigroup $(Q_t)$ and its derivatives imply 
estimates of the semigroup $(P_t)$.

%%%%%
\subsection{Semigroup of subordinate Brownian motion}\label{subsec:subord-BM}
%%%%%

Let $S=(S_t,t\ge 0)$ be a 
subordinator 
that is a nonnegative valued increasing L\'evy process starting at zero. It is characterized by its Laplace 
exponent 
$\phi$ via
\begin{align*}
\bE[\exp(-\lambda S_t)] =e^{-t\phi(\lambda)},\quad t\ge 0 , \lambda>0.
\end{align*}
The Laplace exponent $\phi$ can be written in the form
$$ \phi(\lambda) = b\lambda + \int_{(0,\infty)} (1-e^{-\lambda t}) \mu(dt),$$
where $b\ge 0$ and $\mu$ is a measure on $(0,\infty)$ satisfying 
$\int_{(0,\infty)} 
(1\wedge t) \mu(dt) <\infty$, called the L\'evy measure. Here, $b$ and $\mu(A)$ describes the drift of $S_t$ and the intensity of its jumps of size $A$. In this paper we assume that $b=0, \phi(1)=1$ and $\lim_{\lambda \to \infty} 
\phi(\lambda)=\infty$. Thus $$ \phi(\lambda) = \int_{(0,\infty)} (1-e^{-\lambda 
t})\mu(dt).$$
 
Let $W=(W_t: t \ge 0)$ be the $d$-dimensional Brownian motion with the 
transition density 
$(4\pi t)^{-d/2}\exp{(-|x|^2/(4t))}$ independent of $S$. Define a subordinate 
Brownian 
motion $X=(X_t:t\ge 0)$ by $X_t= W_{S_t}$. The characteristic function of $X$ 
is given by
$$ \bE \big[e^{i\xi \cdot X_t}\big] = e^{-t\phi(|\xi|^2)}$$
and $X$ has the transition density
\begin{align*}
q_d(t,x) = \frac{1}{(2\pi)^d} \int_{\bR^d} e^{i\xi \cdot x}e^{-t\phi(|\xi|^2)} 
d\xi.
\end{align*}
Furthermore, if we denote the distribution of $S_t$ by $\eta_t(dr) =\bP(S_t\in 
dr)$ then 
$q_d(t,x)$ is the same as
\begin{align*}
\int_{(0,\infty)} (4\pi s)^{-d/2} \exp\left(-\frac{|x|^2}{4s}\right) 
\eta_t(ds).
\end{align*}
Thus $q_d(t,x)$ is smooth in $x$. Moreover, its L\'evy measure has a rotationally  symmetric density $J(x)=j(|x|)$ with respect to the Lebesgue measure given by 
$$ j(r) = \int_0^\infty (4\pi t)^{-d/2} e^{-r^2/(4t)} \mu(dt) \,.$$
Note that $J$ is the same function as in \eqref{eq:generator-sub-BM}.

In order to obtain the necessary estimates on the semigroup of subordinate 
Brownian motion we make use of estimates on the transition density and its 
derivatives. In \cite{KKK13} the authors obtain upper bounds of spatial 
derivatives of $q_d(t,x)$ when $\phi$ has a certain scaling condition. For our purposes a weaker version than \cite[Lemma 4.1]{KKK13} is sufficient. We formulate this result without a proof.

\begin{theorem}\label{t:kernel_est}
Suppose $\phi$ satisfies condition \eqref{eq:scaling-cond-phi}. There exists a 
constant $C \geq 1$ such that the following inequalities hold:
\begin{align}
&q_d(t,x) \asymp \left(\phi^{-1}(t^{-1})^{d/2} \wedge 
\frac{t\phi(|x|^{-2})}{|x|^d}\right),\label{e:com_q}\\
&\sum_{|\gamma|= k} \left|D^{\gamma} q_d(t,x) \right| \le C 
\phi^{-1}(t^{-1})^{k/2}\left( 
\phi^{-1}(t^{-1})^{d/2} \wedge 
\frac{t\phi(|x|^{-2})}{|x|^d}\right)\,, \label{e:upperbound_Dq}
\end{align}
for every $k\in\N$ and for all $(t,x) \in(0,\infty)\times \bR^d$.
\end{theorem}

\begin{corollary}\label{cor:bounds-Qt}
Suppose $\phi$ satisfies condition \eqref{eq:scaling-cond-phi} and $k \in 
\N$. There exists 
a constant $C$ depending only on $k, a_1, a_2, \delta_1, \delta_2$ and $d$ such 
that for every multi-index $\gamma$ with $|\gamma| = k$, and every bounded 
function $f$
\begin{align*}
\left| D^\gamma Q_t f(x) \right| \le C \phi^{-1}(t^{-1})^{k/2} \| f\|_{C^0}, 
\quad 
t\in(0,\infty), x\in \bR^d \,.
\end{align*}
\end{corollary}

\pf Comparability of the heat kernel $q_d(t,x)$ and 
$\phi^{-1}(t^{-1})^{d/2} \wedge \left(t\phi(|x|^{-2})|x|^{-d}\right)$ from  
\eqref{e:com_q} and estimate \eqref{e:upperbound_Dq} 
imply
\begin{align*}
\left|D^\gamma Q_tf(x)\right| \le \int_{\bR^d} c_1\phi^{-1}(t^{-1})^{k/2}q_d(t,x-y)f(y) dy 
\le 
c_1\phi^{-1}(t^{-1})^{k/2} \|f\|_{C^0}
\end{align*}
for every multi-index $\gamma$ with $|\gamma| = k$.
\qed 

%%%%%
\subsection{Proof of \autoref{theo:mainL_0-intro}}
%%%%%

The aim of this subsection is to prove \autoref{theo:mainL_0-intro}. Recall 
that $(P_t)$ is the semigroup corresponding to the infinitesimal generator $A$. Let $C_0\ge 1$ be the constant that 
ensures 
\[ C_0^{-1}|h|^{-d}\vphi(|h|)^{-1} \le J(h) \le C_0 |h|^{-d}\vphi(|h|)^{-1} 
\qquad (h \in \R^d\setminus\{0\}) \,. \]
An immediate consequence is that
\[ (\Lambda_2 C_0)^{-1} \; \frac{a_0(h)}{|h|^{d}\vphi(|h|)} \le J(h) \le 
\Lambda_1^{-1} C_0 \; \frac{a_0(h)}{|h|^{d}\vphi(|h|)}  \qquad (h \in \R^d\setminus\{0\}) \,. \]
The derivative estimates of $(Q_t)$ from \autoref{cor:bounds-Qt} imply estimates of 
$(P_t)$ as the next result shows.

\begin{theorem}\label{t:DP_tf}
If $f\in C^0$, then $P_tf$ is $C^{\infty}(\bR^d)$ for $t>0$ and for each 
multi-index 
$\gamma$ with $|\gamma| = k$, there exists $C>0$ (depending on $k$) such that
\begin{align*}
\left| D^\gamma P_tf(x) \right| \le C \vphi^{-1}(t)^{-k} \| f\|_{C^0}.
\end{align*}
\end{theorem}
\pf We define $$ \cL_1 f(x) = \int_{\R^d} (f(x+h)-f(x)-\nabla f(x)\cdot h \1_{|h|\le 1}) (\Lambda_2C_0)^{-1}J(h) \d h $$
and $\cL_2 f(x) = \cL_0f(x)-\cL_1f(x)$ for every $f\in C^2_b(\R^d)$. Let $Q^1_t$ and $Q^2_t$ be the semigroups whose infinitesimal generators are $\cL_1$ and $\cL_2$ respectively. Then $P_t = Q^1_t Q^2_t$. Since $Q^1_t$ is the semigroup of the deterministic time changed process considered in \autoref{t:kernel_est}, we can apply it to $Q^1_t$. Using the contraction property of $Q^2_t$, we get
\begin{align*}
|D^\gamma P_tf(x)| \le c_1\vphi^{-1}(t)^{-k}\|Q^2_tf\|_{C^0} \le c_2\vphi^{-1}(t)^{-k} \|f\|_{C^0}.
\end{align*}
%
%Let $\cL_1$ be the generator 
%of the semigroup corresponding to the L\'evy measure $(\Lambda_1 C_0)^{-1} J(h) \d h$ and let $\cL_2 = \cL_0-\cL_1$. Let $Q^1_t$ and $Q^2_t$ be the semigroups 
%for 
%the L\'evy 
%processes with generators $\cL_1, \cL_2$, respectively, and let $X^1, X^2$ be the 
%corresponding 
%L\'evy processes. If we take $X^1$ independent of $X^2$, then $X^1+X^2$ has the 
%law of 
%the L\'evy process corresponding to the generator $\cL_0$. Therefore $P_t= 
%Q^2_tQ^1_t$. 
%We know that $Q^1_tf$ satisfies the desired estimate by \eqref{e:deriv_Qt} and 
%the fact 
%that the process associated with $\cL_1$ is a deterministic time change of the 
%process 
%considered in \autoref{t:kernel_est}. By translation invariance, $Q^2_t$ 
%commutes with 
%differentiation. Therefore $P_tf = Q^2_tQ^1_tf$ also satisfies the desired 
%estimate, since
%\begin{align}
%\| D^\beta P_t f \|_{C^0} = \| Q^2_t D^\beta Q^1_tf \|_{C^0} \le \| D^\beta 
%Q^1_tf\|_{C^0} \le c \vphi^{-1}(t)^{-k}\| f\|_{C^0}.
%\end{align}
\qed

Recall that we denote the interval of scaling orders of $\psi$ by 
$I_\psi=[m_\psi, M_\psi]$. We  
define the 
potential operator as 
\begin{align*}
R f(x) = \int_0^\infty P_tf(x) dt
\end{align*}
when the function $t \mapsto P_tf(x)$ is integrable. We want to prove that $R$ 
takes 
functions in $C^\psi$ into functions in $C^{\vphi\psi}$, provided that both $I_\psi$ 
and 
$I_{\vphi\psi}$ contain no integer and $Rf$ is bounded.

\begin{lemma}\label{lem:f_eps}
Let $\rho$ be a non-negative $C^\infty$ symmetric function with its support in 
$B(0,1)$ such 
that $\int_{\R^d} \rho(x) dx = 1$, and let $\rho_\eps (x) = \eps^{-d}\rho(x/\eps)$. 
Define 
$f_\eps = f \ast \rho_\eps$. Then for a function $\psi$ with $I_\psi\subset(0,1)$ there 
exists a 
constant $C>0$ such that 
\begin{align}
\| f- f_\eps \|_{C^0} &\le C \| f\|_{C^{\psi}} \psi(\eps),\label{e:f-f_eps}\\
\| D^k f_\eps\|_{C^0} &\le C \|f\|_{C^{\psi}} 
\psi(\eps)\eps^{-k}, \quad k \ge 1.\label{e:Df_eps}
%\| D_{ij} f_\eps\|_{C^0} & \le C \| f\|_{C^{\psi}} \psi(\eps) 
%\eps^{-2}.\label{e:DDf_eps}
\end{align}
\end{lemma}
\pf \eqref{e:f-f_eps} follows from
\begin{align*}
|f(x)-f_\eps(x)| \le  \|f\|_{C^{\psi}} \int_{\R^d} \psi(\eps |y|)\rho(y) dy \le c_1\|f\|_{C^\psi} \psi(\eps).
\end{align*}
In the last inequality we used the fact that $r^{-m_\psi/2}\psi(r)$ is almost increasing. 
Using $\int_{\R^d}D^\gamma \rho (y) dy = 0$ for $|\gamma| \ge 1$,  we can get \eqref{e:Df_eps} from 
\begin{align*}
|D^\gamma f_\eps(x) | &= \left|\int_{\R^d} (f(x-y) - f(x) )D^\gamma\rho_\eps(y) \d y \right|\\
&\le \eps^{-k}\int_{\R^d} |f(x-\eps y)-f(x)| |D^\gamma \rho(y)| \d y \\
&\le c_2\|f\|_{C^\psi} \frac{\psi(\eps)}{\eps^{k}} \int_{\R^d} |D^\gamma \rho(y)| \d y
\end{align*}
for every $x\in \R^d$ and $|\gamma| = k$.
\qed

\begin{proposition}\label{prop:main-proof_prop1}
Suppose $I_\psi\subset(0,1)$ and   
$I_{\vphi\psi}\subset (0,1)\cup (1,2)$. If $f\in C^{\psi}$ and  $\|Rf\|_{C^0}<\infty$, then $Rf\in C^{\vphi\psi}$ and there exists $C$ not depending on $f$ such that $$\|Rf\|_{C^{\vphi\psi}} \le C (\|f\|_{C^{\psi}} + 
\|Rf\|_{C^0}).$$
\end{proposition}
\pf By \autoref{p:equi_norms1} it is enough to show that $[[Rf]]_{C^{\vphi\psi}}$ is 
bounded by $\|f\|_{C^\psi} + \|Rf\|_{C^0}$. Since $|\Delta_h^2(Rf)(x)| \le 4\|Rf\|_{C^0}$ for any $x\in\R^d$, we may assume $|h|\le 1$. 
First we show 
\begin{align}\label{e:Vf}
|\Delta^2_h(P_sf)(x)| \le c_1 |h|^2\|f\|_{C^{\psi}} 
\frac{\psi(\vphi^{-1}(s))}{\vphi^{-1}(s)^2}.
\end{align}
By the Taylor's theorem, \autoref{t:DP_tf}, and \eqref{e:f-f_eps}, 
\begin{align}
|\Delta^2_h(P_s(f-f_\eps))(x) | &
\le |h|^2\| D^2P_s(f-f_\eps)\|_{C^0}\nn\\
&\le \frac{c_2}{\vphi^{-1}(s)^2}|h|^2\|f-f_\eps\|_{C^0}\nn\\
&\le \frac{c_3}{\vphi^{-1}(s)^2}|h|^2\psi(\eps)\| 
f\|_{C^{\psi}}.\label{e:V(f-f_eps)}
\end{align}
Since $\Delta^2_h$ and $P_s$ commute and $P_s$ is a contraction semigroup,   \eqref{e:Df_eps} implies
\begin{align}
|\Delta^2_h(P_sf_\eps)(x)| &=|P_s(\Delta^2_hf_\eps)(x)| \le \| 
\Delta_h^2f_\eps\|_{C^0} \le c_4|h|^2\frac{\psi(\eps)}{\eps^2} 
\|f\|_{C^\psi}.\label{e:Vf_eps}
\end{align}
Letting $\eps =\vphi^{-1}(s)$ and combining with \eqref{e:V(f-f_eps)}, we 
obtain 
\eqref{e:Vf}.

Let $\sigma$ be a small number such that $M_\vphi +M_\psi+2\sigma <2$. Using 
\eqref{e:Vf} and noting $M_\vphi+M_\psi <2$, we have that for $|h|<1$,
\begin{align}
\quad \int_{\vphi(|h|)}^1 |\Delta^2_h(P_sf)(x)| ds &\le 
c_4|h|^2\|f\|_{C^\psi} 
\int_{\vphi(|h|)}^{1} \frac{\psi(\vphi^{-1}(s))}{\vphi^{-1}(s)^2} ds \nn\\
&\le c_5|h|^2\|f\|_{C^{\psi}}\frac{\psi(|h|)}{|h|^2} \int_{\vphi(|h|)}^1 
\left(\frac{|h|}{\vphi^{-1}(s)}\right)^{2- M_\psi-\sigma}ds\nn\\
&\le 
c_6\|f\|_{C^\psi}\psi(|h|)\int_{\vphi(|h|)}
^1 
\left(\frac{\vphi(|h|)}{s}\right)^{(2-M_\psi-\sigma)/(M_\vphi+\sigma)} ds\nn\\
&\le c_7
\vphi(|h|)\psi(|h|)\|f\|_{C^\psi}\label{e:int_Vf_1}
\end{align}

Also the H\"older continuity of $f$ gives
\begin{align*}
 |\Delta^2_h(P_sf)(x)| =|P_s(\Delta^2_hf)(x)| \le 
\|\Delta^2_hf\|_{C^0} \le 
c_8\|f\|_{C^{\psi}}\psi(|h|),
\end{align*}
and thus
\begin{align}\label{e:int_Vf_0}
\int_0^{\vphi(|h|)} |\Delta^2_h(P_sf)(x)| ds \le c_8\|f\|_{C^\psi}\vphi(|h|)\psi(|h|).
\end{align}
Since $|\Delta_h^2 (P_sf)(x)| \le \|D^2P_sf\|_{C^0} |h|^2 \le c_9\vphi^{-1}(s)^{-2}\|f\|_{C^0}$ and $\int_1^\infty \vphi^{-1}(s)^{-2} ds <\infty $, 
\begin{align}\label{e:int_Vf_infty}
\int_1^\infty \left|\Delta_h^2 (P_s f)(x) \right| ds  &\le c_9|h|^2\|f\|_{C^0} \int_1^\infty \vphi^{-1}(s)^{-2} ds  \le c_{10}\vphi(|h|)\psi(|h|)\|f\|_{C^0}.
\end{align}

Adding \eqref{e:int_Vf_1}, \eqref{e:int_Vf_0}, and \eqref{e:int_Vf_infty} we conclude
$$ |\Delta^2_h(Rf)(x)| \le c_{11}\|f\|_{C^\psi} \vphi(|h|)\psi(|h|).$$\qed

Finally we consider the case when $I_{\vphi\psi}\subset (2,3)$.
\begin{proposition}\label{prop:main-proof_prop2}
Suppose $I_\psi\subset(0,1)$, and  
$I_{\vphi\psi}\subset (2,3)$. If $f\in C^{\psi}$ and  $\|Rf\|_{C^0}<\infty$, then $Rf\in C^{\vphi\psi}$ and there exists $C$ not depending on $f$ such that $$\|Rf\|_{C^{\vphi\psi}} \le C (\|f\|_{C^{\psi}} + \|Rf\|_{C^0}).$$
\end{proposition}

\pf Since $I_\psi=[m_\psi, M_\psi]\subset(0,1)$ and 
$I_{\vphi\psi}=[m_\vphi+m_\psi, 
M_\vphi+M_\psi]\subset (2,3)$, necessarily $m_\vphi>1$. Define $\bar{\vphi}(r) =r^{-1}\vphi(r)$ then $I_{\bar{\vphi}\psi} \subset(1,2)$. In view of 
\autoref{p:equi_norms1} it suffices to show 
\begin{align}\label{e:DRf}
[[ D Rf ]]_{C^{\bar{\vphi}\psi}} \le c_1 \| f\|_{C^\psi}.
\end{align}
Fix $i$ and let $Q_sf(x) = D_i (P_sf)(x)$. From \autoref{t:DP_tf} we have
\begin{align*}
\| D^2 Q_s f\|_{C^0} \le c_2 \vphi^{-1}(s)^{-3} \| f\|_{C^0}.
\end{align*}
Note that $Q_s$ is translation 
invariant. As the proof of \autoref{prop:main-proof_prop1}
 we assume $|h|\le 1$. Analogously to \eqref{e:V(f-f_eps)} and \eqref{e:Vf_eps},
\begin{align*}
|\Delta^2_h(Q_s(f-f_\eps))(x) | 
%&\le |h|^2\| D^2Q_s(f-f_\eps)\|_{C^0}\\
%&\le \frac{c_3}{\vphi^{-1}(s)^3} |h|^2\|f-f_\eps\|_{C^0}\\
%&
\le \frac{c_4}{\vphi^{-1}(s)^3} |h|^2\psi(\eps)\|f\|_{C^{\psi}}
\end{align*}
and
\begin{align*}
|\Delta^2_h(Q_sf_\eps)(x)| &\le |h|^2 \| D^2Q_s f_\eps\|_{C^0} \le c_5|h|^2\|D^3 f_\eps \|_{C^0} \le 
c_6|h|^2
\frac{\psi(\eps)}{\eps^3} \| f\|_{C^{\psi}}.
\end{align*}
Taking $\eps = \vphi^{-1}(s)$ we obtain
$$ |\Delta^2_h(Q_sf)(x)| \le c_6|h|^2 \frac{\psi(\vphi^{-1}(s))}{\vphi^{-1}(s)^3} 
\|f\|_{C^{\psi}}.$$
Integrating the right hand side with respect to $s$ over the interval $[\vphi(|h|),1)$ yields 
$c_7\bar{\vphi}(|h|)\psi(|h|)\|f\|_{C^{\psi}}$.

On the other hand, 
\begin{align*}
|\Delta^2_h(Q_sf)(x)| \le c_8 
\vphi^{-1}(s)^{-1}\| 
\Delta_h^2f\|_{C^0} \le 2c_8\vphi^{-1}(s)^{-1}\psi(|h|)\|f\|_{C^\psi}.
\end{align*}
and integrating this bound over the interval $(0,\vphi(|h|))$ yields 
$c_9\bar{\vphi}(|h|)\psi(|h|)\|f\|_{C^{\psi}}$; we use  
$m_\vphi>1$ here.
Since $|\Delta_h^2 Q_sf(x)| \le c_{10}|h|^2\|f\|_{C^0} \vphi^{-1}(s)^{-3} $ and $\int_1^\infty \vphi^{-1}(s)^{-3} ds <\infty$,
\begin{align*}
\int_1^\infty|\Delta_h^2(Q_sf)(x)| ds \le c_{10}|h|^2\|f\|_{C^0} \int_1^\infty \vphi^{-1}(s)^{-3} ds \le c_{11}\bar{\vphi}(|h|)\psi(|h|)\|f\|_{C^0}.
\end{align*} 
Therefore
$$ |\Delta^2_h(D_i Rf)(x)| \le c_{12}\bar{\vphi}(|h|)\psi(|h|)
\|f\|_{C^{\psi}},$$
which yields \eqref{e:DRf}.\qed

Now, the proof of \autoref{theo:mainL_0-intro} follows by the preceding 
propositions.

\pf [Proof of \autoref{theo:mainL_0-intro}] According to the assumptions of 
\autoref{theo:mainL_0-intro} the function $u$ is an element of 
$C^{\vphi \psi}(\R^d)\cap C_\infty(\R^d)$. Without loss of generality we may assume that $u$ 
belongs to $C^2_\infty(\R^d) \cap C^{\vphi \psi}(\R^d)$ and thus to the domain 
$\cD(A)$ of the infinitesimal generator $(A, \cD(A))$ of 
the semigroup $(P_t)$. Because we may convolve $u$ with a mollifier 
$\rho_\eps$ like in \autoref{lem:f_eps}. Then $u_\eps = u * 
\rho_\eps$ is a smooth function 
vanishing at infinity and satisfies the equation $\cL_0 
u_\eps =  f * \rho_\eps$. We would then obtain the estimate claimed in 
\autoref{theo:mainL_0-intro} for $u_\eps$. Since the three norms in this estimate 
converge for $\eps \to 0$, the desired estimate for $u$ would 
follow. 

Recall that the 
infinitesimal generator $(A, \cD(A))$ of 
the semigroup $(P_t)$ satisfies $A v = \cL_0 v$ for every $v \in 
C^2_\infty(\R^d)$ 
and $\cL_0$ as in \eqref{def:operator_const_coeff}. Denote by $(R, \cD(R))$ be 
the potential operator of $(P_t)$, i.e.,
\[ Rf = \lim_{t \to \infty} \int_0^t P_s f \d s \,. \]
Note that in general the potential operator is not 
identical with the zero-resolvent operator $(R_0, \cD(R_0))$. However, the property that
\[ \|P_t v\|_{C^0} \to 0 \text{ as }  t \to \infty \quad  \text{ for every } v \in C_\infty  \]
from the translation invariance implies that $(R, \cD(R))$ is 
densely defined and $R=R_0 = -A^{-1}$ 
\cite[Proposition 11.9]{BeFo75}. Hence $u =- Rf$ and we can apply 
\autoref{prop:main-proof_prop1} and \autoref{prop:main-proof_prop2} from above.
\qed

%%%%%%%%%%%%%%%%%%%%%%%%%%%%%%%%%%%%%%%%%%%%%%%%%%%%%%%%%%%%%%%%%%%%%%%%%%%%
%%%%%%%%%%%%%%%%%%%%%%%%%%%%%%%%%%%%%%%%%%%%%%%%%%%%%%%%%%%%%%%%%%%%%%%%%%%%
\section{Proof of \autoref{theo:main}}\label{sec:main-proof}
%%%%%%%%%%%%%%%%%%%%%%%%%%%%%%%%%%%%%%%%%%%%%%%%%%%%%%%%%%%%%%%%%%%%%%%%%%%%
%%%%%%%%%%%%%%%%%%%%%%%%%%%%%%%%%%%%%%%%%%%%%%%%%%%%%%%%%%%%%%%%%%%%%%%%%%%%

The aim of this section is to prove \autoref{theo:main} using 
\autoref{theo:mainL_0-intro} and a well-known perturbation technique. Let us 
first establish some auxiliary results.

We show that \eqref{eq:scaling-cond-vphi} implies 
\begin{align}\label{ass:levymeasure_outside}
\int_r^\infty \frac{\d s }{s \vphi(s)} \le \frac{C}{\vphi(r)}  \qquad (r > 
0)\,,
\end{align}
where $C$ is some positive constant. 
The second inequality in 
\eqref{eq:scaling-cond-vphi} with $\lambda = s/r$ implies 
\[ \vphi(s) \geq  a_1 (s/r)^{2 
\delta_1}  \varphi(r) \,. \]

The above observation \eqref{ass:levymeasure_outside} now follows from
\begin{align*}
\int_r^\infty \frac{\d s }{s \vphi(s)} 
\le a_1^{-1} \frac{r^{2 \delta_1}}{\vphi(r)} \int_r^\infty \frac{\d s 
}{s^{1+2\delta_1} }  = \frac{1}{2a_1\delta_1 \vphi(r)}  \,. 
\end{align*}

Let $B(x,r)$ denote the ball of radius $r$ centered at $x$. Let $\bar{\eta} \in 
C^\infty_c(\R^d)$ be 
a cut-off 
function which equals $1$ on $B(0,1)$, equals $0$ on $B(0,2)^c$ and 
satisfies $\bar{\eta} \in [0,1]$. Let $\eta_{r,x_0}(x) = 
\bar{\eta}((x-x_0)/r)$. When there is no ambiguity we write $\eta$ instead of 
$\eta_{r,x_0}$. 

\begin{proposition}\label{p:u_eta}
Assume $I_\psi\subset(0,1)$, $I_{\vphi\psi}\subset(0,1)\cup (1,2)\cup(2,3)$,  $u \in C^{\vphi\psi}(\R^d)$ and  $f \in C^{\psi}(\R^d)$. Suppose that 
for each  $\eps>0$ there exist $r>0$ and $c_1\geq 1$ depending on $\eps$ such that 
\begin{align}\label{e:u_eta}
\|u\eta_{r,x_0}\|_{C^{\vphi\psi}} \le c_1(\|f\|_{C^\psi} + \| u \|_{C^0} )  + 
\eps \| 
u\|_{C^{\vphi\psi}}
\end{align}
for all 
$x_0\in\bR^d$. Then there exists a constant $C$ such that
\begin{align*}
\|u\|_{C^{\vphi\psi}} \le C(\|f\|_{C^\psi} + \|u\|_{C^0}).
\end{align*}
\end{proposition}
\pf First, we consider the case $I_{\vphi\psi}\subset(0,1)\cup(1,2)$. Set  
$\eps=1/2$ and choose $r$ and $c_1$ satisfying \eqref{e:u_eta}. For any $x_0\in 
\bR^d$, if  
$|h|<r$ then $\Delta_h^2u(x_0) = \Delta_h^2(u\eta_{r,x_0})(x_0)$, and \eqref{e:u_eta} yields
\begin{align*}
|\Delta^2_hu(x_0)| & \le 
\|u\eta_{r,x_0}\|_{C^{\vphi\psi}}\vphi(|h|)\psi(|h|) \le \left(c_1\|f\|_{C^\psi} + c_1\| u\|_{C^0}   + \frac{1}{2} \| u\|_{C^{\vphi\psi}}\right) \vphi(|h|)\psi(|h|)
\end{align*}
On the other hand, if $r\le |h| \le 1$ then the  fact that $\vphi(s)\psi(s)\ge c_2$ for any $s\in[r, 1]$ yields
\begin{align*}
|\Delta^2_hu(x_0)| \le 4c_2^{-1} \|u\|_{C^0} 
\vphi(|h|)\psi(|h|).
\end{align*}
Combining the above two inequalities we obtain
\begin{align*}
\sup_{x\in\bR^d}\sup_{0<|h|\le 1} \frac{|\Delta_h^2u(x)|}{\vphi(|h|)\psi(|h|)} 
\le 
c_1\|f\|_{C^\psi} + c_1(1+4c_2^{-1})\|u\|_{C^0} + \frac{1}{2} \|u\|_{C^{\vphi\psi}}
\end{align*}
Therefore we obtain
$$ \|u\|_{C^{\vphi\psi}} \le c_3\left(\|f\|_{C^\psi} + \|u\|_{C^0} \right).$$

Now we consider the case when $I_{\vphi\psi}\subset(2,3)$. Let $\bar{\vphi}(r) = r^{-1}\vphi(r)$.  By the definition of 
$C^{\vphi\psi}$ and \autoref{p:equi_norms1} it is enough to show that $$\|Du\|_{C^0} + [[Du]]_{C^{\bar{\vphi}\psi}}\le c_4(\|f\|_{C^\psi}+\|u\|_{C^0}).$$ 
For $\eps = 1/4$ choose $r$ satisfying \eqref{e:u_eta}. We use the same argument above to obtain that if 
$|h|<r$ then 
\begin{align*}
|\Delta^2_h(Du)(x_0)| 
%&= |\Delta^2_h (Dv)(x_0)| \le [[Dv]]_{C^{\bar{\vphi}\psi}} 
%\bar{\vphi}(|h|)\psi(|h|)\\
%&\le \|v\|_{C^{\vphi\psi}} \bar{\vphi}(|h|)\psi(|h|) \\
&\le \left(c_5\|f\|_{C^\psi} + c_5\|u\|_{C^0} + \frac{1}{4}\|u\|_{C^{\vphi\psi}}\right)\bar{\vphi}(|h|)\psi(|h|)
\end{align*}
for any $x_0\in\bR^d$. On the other hand, if $r\le |h|\le 1$ then the  fact that $\bar{\vphi}(s)\psi(s)\ge c_6$ for any $s\in[r,1]$ yields
\begin{align*}
|\Delta^2_h(Du)(x_0)| \le 4c_6^{-1}\|Du\|_{C^0} 
\bar{\vphi}(|h|)\psi(|h|).
\end{align*}
Combining above two inequalities and we get
\begin{align*}
\sup_{x\in\bR^d}\sup_{0<|h|\le 1} \frac{|\Delta_h^2(Du)(x)|}{\bar{\vphi}(|h|)\psi(|h|)} 
\le 
c_7\left(\|f\|_{C^\psi} + \|u\|_{C^0} + \|Du\|_{C^0} \right)+ \frac{1}{4} \|u\|_{C^{\vphi\psi}}
\end{align*}

By \autoref{p:ineq_norms} we have $$ \|Du\|_{C^0} \le c_8\|u\|_{C^0} + 
(4(1+c_7))^{-1}\|u\|_{C^{\vphi\psi}},$$
which implies 
\begin{align*}
 \|Du\|_{C^0} + 
\sup_{x\in\bR^d}\sup_{0<|h|\le 1}\frac{|\Delta^2_h(Du)(x)|}{\bar{\vphi}(|h|)\psi(|h|)} &\le c_9\left(\|f\|_{C^\psi} + \|u\|_{C^0}\right) + \frac{1}{2}\|u\|_{C^{\vphi\psi}}.
\end{align*}
Therefore we obtain the desired estimate
\begin{align*}
\|u\|_{C^{\vphi\psi}} \le c_{10}(\|f\|_{C^\psi} + \|u\|_{C^0}).
\end{align*}
\qed

Before proving the main theorem, we give an auxiliary inequality, which we will often apply.
\begin{lemma}\label{l:aux_int}
Let $\Psi :(0,1] \to (0,\infty)$ be a function with $M_\vphi < m_\Psi$. There exists a constant $C>0$ such that for any $0<r\le 1$,
\begin{align}\label{e:aux_int}
\int_{\R^d} \frac{\Psi(|h|\wedge r)}{|h|^d\vphi(|h|)} \d h \le C\frac{\Psi(r)}{\vphi(r)}.
\end{align}
\end{lemma}
\pf Let $\sigma = (m_\Psi-M_\vphi)/3>0$. By the definition of $I_\vphi$ and $I_\Psi$, there exists a constant $c_1>0$ such that for $0<|h|< r$,
\begin{align*}
\frac{\Psi(|h|)}{\vphi(|h|)} \le c_1 \frac{(|h|/r)^{m_\Psi-\sigma}\Psi(r)}{(|h|/r)^{M_\vphi+\sigma}\vphi(r)} \le c_1 (|h|/r)^\sigma\frac{\Psi(r)}{\vphi(r)}.
\end{align*}
It follows that 
\begin{align*}
\int_{|h|<r} \frac{\Psi(|h|)}{|h|^d\vphi(|h|)} \d h \le c_1\frac{\Psi(r)}{\vphi(r)} \int_{|h|<r} \frac{(|h|/r)^\sigma}{|h|^d} \d h \le c_2\frac{\Psi(r)}{\vphi(r)}. 
\end{align*}
Combining this with \eqref{ass:levymeasure_outside}, we get the result. \qed

Let $H$ be a function defined by 
\begin{align*}
H(x) = \int_{\R^d} (u(x+h)-u(x))(\eta(x+h)-\eta(x))\frac{a(x,h)}{|h|^d\vphi(|h|)} \d h.
\end{align*}

\begin{lemma}\label{l:H}
Let $\eps>0$ be a small constant. Assume that condition \eqref{ass:I_contain_no_integer} is satisfied. If $u\in C^{\vphi\psi}$, then  there exists a constant $C=C(r,\eps)>0$ such that
\begin{align*}
\| H\|_{C^\psi} \le C\| u\|_{C^0} + \eps \|u\|_{C^{\vphi\psi}}.
\end{align*}
\end{lemma}
\pf Observe that 
\begin{align}\label{e:pf_Dhu}
|\Delta_h u(x) | \le 
\begin{cases}
2 \|u\|_{C^0},\quad &\text{ if } I_{\vphi\psi}\subset (0,1),\\
2(\|u\|_{C^0}+\|Du\|_{C^0})(|h|\wedge1), &\text{ if } I_{\vphi\psi}\subset(1,2)\cup(2,3),
\end{cases}
\end{align}
and 
\begin{align}\label{e:pf_Dh_eta}
|\Delta_h \eta(x) | \le 
\begin{cases}
2\|\bar{\eta}\|_{C^{\vphi\psi}}(\vphi\psi)\left(\frac{|h|}{r} \wedge 1\right), \quad  & \text{ if } I_{\vphi\psi} \subset(0,1), \\
2(\|\bar{\eta}\|_{C^0} + \|D\bar{\eta}\|_{C^0}) \left(\frac{|h|}{r} \wedge 1\right), &  \text{ if } I_{\vphi\psi} \subset(1,2)\cup(2,3).  
\end{cases}
\end{align}
When $I_{\vphi\psi}\subset(0,1)$, using \eqref{e:pf_Dhu} and \eqref{e:pf_Dh_eta} we get
\begin{align*}
|H(x)| &\le \int_{\R^d} |\Delta_h u(x)| |\Delta_h \eta(x)| \frac{|a(x,h)|}{|h|^d\vphi(|h|)} \d h \\
&\le 4\|u\|_{C^0}\|\bar{\eta}\|_{C^{\vphi\psi}}\int_{\R^d} (\vphi\psi)\left(\frac{|h|}{r}\wedge 1\right) \frac{\Lambda_2 }{|h|^d\vphi(|h|)} \d h \\
&\le \frac{c_1}{\vphi(r)} \|u\|_{C^0}.
\end{align*}
We used \eqref{e:aux_int} with $\Psi(t) = (\vphi\psi)(t/r)$ in the last inequality. When $I_{\vphi\psi} \subset(1,2)\cup(2,3)$,  by \eqref{e:pf_Dhu}, \eqref{e:pf_Dh_eta}, and \eqref{e:aux_int} with $\Psi(t) = (t/r)^2$, we get 
\begin{align*}
|H(x)| \le \frac{c_2}{\vphi(r)} (\|u\|_{C^0} + \|Du\|_{C^0}). 
\end{align*}
By \autoref{p:ineq_norms}, we can find a constant $c_3= c_3(r,\eps)>0$ such that
\begin{align*}
\|H\|_{C^0} \le c_3(r,\eps) \|u\|_{C^0} + \eps\|u\|_{C^{\vphi\psi}},
\end{align*}
whenever $I_{\vphi\psi}\subset(0,1)\cup(1,2)\cup(2,3)$. 

Now, we consider $|\Delta_k H(x)|$ in order to estimate the $\psi$-H\"older seminorm of $H$. We may assume $|k|<r$ because $|\Delta_k H(x)|/\psi(|k|) \le c_4\|H\|_{C^0}/\psi(r)$ otherwise. Observe that
\begin{align*}
\Delta_k H(x) &= \int_{\R^d} \Delta_k\Delta_h u(x) \, \Delta_h \eta(x)\,  a(x+k, h)\, \frac{\d h}{|h|^d\vphi(|h|)}  \\
&\quad +\int_{\R^d} \Delta_h u(x)\, \Delta_k\Delta_h\eta(x) \, a(x+k,h)\, \frac{\d h }{|h|^d\vphi(|h|)}  \\
&\quad +\int_{\R^d} \Delta_h u(x) \, \Delta_h \eta(x) \, \Delta_k (a(\cdot,h))(x)\, \frac{\d h}{|h|^d\vphi(|h|)}  \\
&=: I_1+I_2+I_3.
\end{align*}
Thus it suffices to show that 
\begin{align}\label{e:pf_H_aux}
 |I_1|+|I_2|+|I_3| \le \psi(|k|)( c_5 \|u\|_{C^0} +\eps \|u\|_{C^{\vphi\psi}})
 \end{align}
for some constant $c_5=c_5(r,\eps)>0$.

For $|I_1|$, it follows from the facts $u\in C^{\vphi\psi}$ and  $\eta \in C^\infty$ that
\begin{align*}
|\Delta_k\Delta_hu(x) \, \Delta_h \eta(x) | \le c_6\cdot \begin{cases} 
 \|u\|_{C^\psi}\psi(|k|)\left(\frac{|h|}{r}\wedge 1\right), \quad & \text{ if } I_{\vphi\psi}\subset(0,1),\\
 \|u\|_{C^{\sigma}} |k| \, \frac{\sigma(|h|\wedge 1)}{|h|\wedge 1} \left(\frac{|h|}{r}\wedge 1\right), \quad & \text{ if } I_{\vphi\psi}\subset(1,2),\\
 (\|Du\|_{C^0}+\|D^2u\|_{C^0}) |k| (|h|\wedge 1) \left(\frac{|h|}{r}\wedge 1\right), \quad & \text{ if } I_{\vphi\psi}\subset(2,3),
\end{cases}
\end{align*}
for some constant $c_6>0$, where $\sigma$  is a function on $(0,1]$ defined by
$$\sigma(t) = t^{(1\vee M_\vphi + m_\vphi+m_\psi)/2}.$$
Note that the exponent of $\sigma$ is greater than both $1$ and $M_\vphi$. Since $r$ is less than one, we can apply \autoref{l:aux_int}  with $\Psi_1(t) = t/r$, $\Psi_2(t) = \sigma(t/r)$, and $\Psi_3(t) = (t/r)^2$ for each cases. Then we have
\begin{align*}
|I_1| \le \frac{c_7\Lambda}{\vphi(r)} \cdot 
\begin{cases}
\|u\|_{C^\psi} \psi(|k|), \quad &\text{ if } I_{\vphi\psi} \subset(0,1),\\
\|u\|_{C^{\sigma}} |k|, \quad &\text{ if } I_{\vphi\psi} \subset(1,2),\\
(\|Du\|_{C^0}+\|D^2 u\|_{C^0} )|k|, \quad &\text{ if } I_{\vphi\psi} \subset(2,3),
\end{cases}
\end{align*}
for some constant $c_7>0$. By the fact $|k| \le c_8\psi(|k|)$ and \autoref{p:ineq_norms}, there exists a constant $c_9 = c_9(r,\eps)>0$ such that 
\begin{align*}
|I_1| \le \psi(|k|) (c_9(r,\eps) \|u\|_{C^0} + \eps\|u\|_{C^{\vphi\psi}})
\end{align*}
whenever $I_{\vphi\psi}\subset(0,1)\cup(1,2)\cup(2,3)$.

For $|I_2| + |I_3|$, we first consider the case $I_{\vphi\psi}\subset(0,1)$. 
%Note that $I_\vphi \subset(0,1)$ follows obviously. It is trivial that
%\begin{align}\label{e:pf_H_Dhu}
%|\Delta_h u(x) | \le 2\|u\|_{C^0}.
%\end{align}
Since $|\Delta_k\Delta_h \eta (x) | \le c_{10} \left(\frac{|h|}{r}\wedge 1\right) \psi\left(\frac{|k|}{r}\right)$, we have
\begin{align}\label{e:pf_H_I2_int}
|\Delta_h u(x) \Delta_k\Delta_h \eta(x) a(x+k,h) | \le 2c_{10} \Lambda_2 \|u\|_{C^0}  \left(\frac{|h|}{r}\wedge 1\right)\psi\left(\frac{|k|}{r}\right).
\end{align}
Similarly, the facts that $|\Delta_h \eta(x)| \le c_{11} \left(\frac{|h|}{r}\wedge 1\right)$ and \eqref{eq:assum_axh-cont} implies
\begin{align}\label{e:pf_H_I3_int}
|\Delta_h u(x) \Delta_h \eta(x) \Delta_k(a(\cdot,h))(x) | \le 2c_{11}\Lambda_3 \|u\|_{C^0} \left(\frac{|h|}{r}\wedge 1\right) \psi(|k|). 
\end{align} 
By integrating \eqref{e:pf_H_I2_int}, \eqref{e:pf_H_I3_int}, and then  applying \autoref{l:aux_int} with $\Psi(t) = t/r$, we can obtain 
\begin{align*}
|I_2|+|I_3| \le \frac{c_{12}}{\vphi(r)} \|u\|_{C^0} \psi(|k|)
\end{align*}
for some constant $c_{12}>0$. When $I_{\vphi\psi} \subset(1,2)\cup(2,3)$, we just use $$|\Delta_h u(x)| \le 2(\|u\|_{C^0} + \|Du\|_{C^0})(|h|\wedge 1),$$
and apply \autoref{l:aux_int} with $\Psi(t) = (t/r)^2$. Finally, we use \autoref{p:ineq_norms} to obtain \eqref{e:pf_H_aux}
%$$|I_2|+|I_3| \le \psi(|k|) (c_{13} \|u\|_{C^0} + \eps \|u\|_{C^{\vphi\psi}})$$
%for some constant $c_{13}=c_{13}(r,\eps)>0$ 
whenever $I_{\vphi\psi}\subset(0,1)\cup(1,2)\cup(2,3)$.
\qed

Let $x_0$ be a fixed point in  $\R^d$ and $b(x,h) = a(x,h)-a(x_0,h)$. When $M_\vphi<1$,  we define $\cB v(x)$ by
\begin{align*}
\cB v(x) = \int_{\R^d} (v(x+h) - v(x)) \frac{b(x,h)}{|h|^d\vphi(|h|)} \d h.
\end{align*}
When $M_\vphi \ge 1$, we define $\wt{\cB} v(x)$  by adding a gradient term as
\begin{align*}
\wt{\cB} v(x) = \int_{\R^d} (v(x+h) - v(x)-\nabla v(x)\cdot h \1_{\{|h|\le 1\}}) \frac{b(x,h)}{|h|^d\vphi(|h|)} \d h.
\end{align*}
\begin{lemma}\label{l:Bv}
Let $\eps>0$ be given. Assume $M_\vphi<1$, $I_{\vphi\psi}\subset(0,1)\cup(1,2)$, and that condition \eqref{ass:I_contain_no_integer} is satisfied. Then there exists $r=r(\eps)\in(0,1/4)$ such that for every $v\in C^{\vphi\psi}$ with its support in $B(x_0, 2r)$, 
\begin{align*}
\|\cB v \|_{C^\psi} \le C\|v\|_{C^0} + \eps \|v\|_{C^{\vphi\psi}}
\end{align*}
for $C= C(\eps)>0$
\end{lemma}
\pf Let $v$ be a function with its support in $B(x_0,2r)$ for $0<r<1/4$. We first obtain a bound of the $C^0$-norm of $\cB v$. If $x\notin B(x_0,3r)$, then $v(x) = 0$ and $v(x+h) = 0$ for $|h|\le r$. Thus \eqref{ass:levymeasure_outside} yields
\begin{align}\label{e:proof-abs_cBv1}
|\cB v(x) | = \left| \int_{|h|>r} v(x+h)\frac{b(x,h)}{|h|^d\vphi(|h|)} \d 
h\right| 
\le 
\|v\|_{C^0} \int_{|h|>r} \frac{2\Lambda_2 \d h}{|h|^d\vphi(|h|)} \le 
\frac{c_1}{\vphi(r)}\|v\|_{C^0}.
\end{align}
Let us look at the case $x\in B(x_0,3r)$. Note that, in this case,  $|b(x,h)|\le c_2\Lambda_3\psi(r)$ for every $h\in\R^d$. Observe that
\begin{align}\label{e:pf_Dhv}
|\Delta_h v(x) | \le 
\begin{cases}
2\|v\|_{C^{\sigma}} \sigma(|h|\wedge 1), \quad & \text{ if } I_{\vphi\psi}\subset(0,1),\\
2(\|v\|_{C^0} + \|Dv\|_{C^0}) (|h|\wedge 1), \quad &\text{ if } I_{\vphi\psi}\subset(1,2),
\end{cases}
\end{align}
where $\sigma$ is a function on $(0,1]$ defined by 
$$ \sigma(t) = t^{(M_\vphi+m_\vphi+m_\psi)/2}.$$
Note that $M_\vphi < m_{\sigma} = (M_\vphi + m_\vphi+m_\psi)/2$. Applying \autoref{l:aux_int} to its integration with $\Psi(t) = \sigma(t)$ when $I_{\vphi\psi}\subset(0,1)$, and with $\Psi(t) = t$ when $I_{\vphi\psi}\subset(1,2)$, we get
\begin{align}\label{e:pf_Bv_abs}
|\cB v(x) | \le c_3\psi(r) \cdot
\begin{cases}
\|v\|_{C^\sigma} , \quad & \text{ if } I_{\vphi\psi}\subset(0,1),\\
\|v\|_{C^0} + \|Dv\|_{C^0}, \quad &\text{ if } I_{\vphi\psi}\subset(1,2).
\end{cases}
\end{align}
It follows from \eqref{e:proof-abs_cBv1} and \autoref{p:ineq_norms}  that
\begin{align}\label{e:pf_Bv}
 \|\cB v\|_{C^0} \le c_4 \|v\|_{C^0} + (\eps/5) \|v\|_{C^{\vphi\psi}}
 \end{align}
for some constant $c_4 = c_4(r,\eps)>0$. 

In the next step we estimate the $\psi$-H\"older seminorm of $\cB v$. If $r/2 < |k| \le 1$, then \eqref{e:proof-abs_cBv1} and \eqref{e:pf_Bv_abs} yield
\begin{align*}
\frac{|\Delta_k(\cB v) (x) |}{\psi(|k|)}  \le \frac{c_5 \|v\|_{C^0} }{\vphi(r)\psi(r/2)}+ \frac{c_5\psi(r)}{\psi(r/2)} \cdot
\begin{cases}
\|v\|_{C^\sigma}, \quad &\text{ if } I_{\vphi\psi}\subset(0,1), \\
(\|v\|_{C^0} + \|Dv\|_{C^0}), \quad &\text{ if } I_{\vphi\psi}\subset(1,2),
\end{cases}
\end{align*}
By the fact $\psi(r) \le c_6 \psi(r/2)$ and \autoref{p:ineq_norms}, the quotient $\frac{|\Delta_k \cB(x)|}{\psi(|k|)}$ for $r/2<|k|\le 1$ is bounded by 
\begin{align}\label{e:pf_Bv_semi1}
 c_7(r,\eps) \|v\|_{C^0} + (\eps/5) \|v\|_{C^{\vphi\psi}}.
 \end{align}
Now consider the case  $|k|\le r/2$. First suppose $x\notin B(x_0,3r)$. Then 
$v(x+k)=v(x)=0$ and $v(x+k+h)=v(x+h)=0$ for $|h|\le r/2$. Thus, the inequality 
$|\Delta_k v(x+h)| \le \|v\|_{C^\psi}\psi(|k|)$, 
 \eqref{ass:levymeasure_outside} and 
 \eqref{eq:assum_axh-cont} yield
\begin{align}
|\Delta_k(\cB v)(x) |&= \left|\int_{|h|>r/2} 
(v(x+k+h)b(x+k,h)-v(x+h)b(x,h))\frac{\d h}{|h|^d\vphi(|h|)}\right| \nn\\
&\le \int_{|h|>r/2} |\Delta_kv(x+h)|\frac{2\Lambda _2\, \d h}{|h|^d\vphi(|h|)} + 
\int_{|h|>r/2} |v(x+h)|\frac{|\Delta_k(b(\cdot,h))(x)|}{|h|^d\vphi(|h|)}\d h\nn \\
&\le 2(\Lambda_2+\Lambda_3) (\|v\|_{C^\psi} + \|v\|_{C^0} )\psi(|k|)\int_{|h|>r/2} \frac{\d h}{|h|^d\vphi(|h|)} \nn\\
&\le \frac{c_{8}}{\vphi(r/2)} \|v\|_{C^\psi}\psi(|k|) \nn
\end{align}
By \autoref{p:ineq_norms} we get 
\begin{align}\label{e:pf_Bv_semi2}
|\Delta_k(\cB v) (x)| \le (c_9(r,\eps) \|v\|_{C^0} + (\eps/5)\|v\|_{C^{\vphi\psi}})\psi(|k|).
\end{align}
Now suppose $x\in B(x_0,3r)$. We decompose the integral into two parts as follows 
\begin{align*}
\Delta_k(\cB v)(x) & =  \int_{\R^d} 
\Delta_k\Delta_hv(x)\frac{b(x+k,h)}{|h|^d\vphi(|h|)}\d h + \int_{\R^d} 
\Delta_hv(x)\frac{\Delta_k(b(\cdot,h))(x)}{|h|^d\vphi(|h|)} \d h \\
&=: I_4+I_5 \,.
\end{align*}
For $I_4$, we observe that  $x+k\in B(x_0,4r)$ 
and $|b(x+k,h)|\le c_{10} \psi(r)$ since $|k|\le r/2$. When $I_{\vphi\psi} \subset(0,1)$, by the inequality 
\begin{align}\label{e:pf_DkDhv}
|\Delta_k\Delta_h v(x)| \le 2\|v\|_{C^{\vphi\psi}} (\vphi\psi)(|h|\wedge |k|),
\end{align}
and \eqref{e:aux_int} with $\Psi(t) = (\vphi\psi)(t)$, we obtain
\begin{align*}
|I_4| &\le 2c_{10}
\psi(r)\|v\|_{C^{\vphi\psi}}\int_{\R^d}\frac{(\vphi\psi)(|h|\wedge |k|)}{|h|^d\vphi(|h|)}\d h \le c_{11}\psi(r)\|v\|_{C^{\vphi\psi}}\psi(|k|).
\end{align*}
If $I_{\vphi\psi}\subset(1,2)$,  the following inequality is used instead of \eqref{e:pf_DkDhv},
\begin{align*}
|\Delta_k\Delta_h v(x)| \le 2\cdot 
\begin{cases}
\|v\|_{C^{\vphi\psi}} |h|\frac{(\vphi\psi)(|k|)}{|k|}, \quad & \text{ if } |h|\le |k|,\\
\|v\|_{C^{\vphi\psi}} |k|\frac{(\vphi\psi)(|h|)}{|h|}, \quad & \text{ if } |k|< |h|\le 1,\\
\|v\|_{C^\psi} \psi(|k|), \quad & \text{ if }  |h|>1.\\
\end{cases}
\end{align*}
By the integration and \autoref{p:ineq_norms} we obtain
\begin{align}\label{e:pf_I4}
|I_4| \le c_{12}\psi(r)(\|v\|_{C^0}+\|v\|_{C^{\vphi\psi}}) \psi(|k|),
\end{align} 
whenever $I_{\vphi\psi}\subset(0,1)\cup(1,2)$.

Since $|\Delta_k(b(\cdot,h))(x)| \le 2\Lambda \psi(|k|)$ for every $h\in\R^d$, \eqref{e:pf_Dhv} and the continued argument yields 
\begin{align}\label{e:pf_I5}
|I_5| \le (c_{13} \|v\|_{C^0} + (\eps/5)\|v\|_{C^{\vphi\psi}})\psi(|k|),
\end{align}
whenever $I_{\vphi\psi}\subset(0,1)\cup(1,2)$. Combining \eqref{e:pf_Bv}, \eqref{e:pf_Bv_semi1}, \eqref{e:pf_Bv_semi2}, \eqref{e:pf_I4} and \eqref{e:pf_I5}, we have
\begin{align*}
\|\cB v\|_{C^\psi} \le c_{14}(r,\eps) \|v\|_{C^0} + c_{12}\psi(r)\|v\|_{C^{\vphi\psi}} + (4\eps/5)\|v\|_{C^{\vphi\psi}}
\end{align*}
Choosing $r$ such that $c_{12}\psi(r) \le \eps/5$, we can obtain the result.
\qed

\begin{lemma}\label{l:BBv}
Let $\eps >0$ be given. Assume $M_\vphi\ge 1$, $I_{\vphi\psi}\subset(1,2)\cup(2,3)$, and the condition \eqref{ass:I_contain_no_integer} is satisfied. In the case $1 \in I_\vphi$ we 
further assume $a(x,h)=a(x,-h)$ for all $x,h\in\R^d$. Then there exists $r=r(\eps)\in(0,1/4)$ such that for every $v\in C^{\vphi\psi}$ with its support in $B(x_0, 2r)$, 
\begin{align*}
\|\wt{\cB} v \|_{C^\psi} \le C\|v\|_{C^0} + \eps \|v\|_{C^{\vphi\psi}}
\end{align*}
for $C= C(\eps)>0$.
\end{lemma}
\pf For the $C^0$-norm of $\wt{\cB} v$, the only difference from the proof of \autoref{l:Bv} is that we replace \eqref{e:pf_Dhv} with 
\begin{align}\label{e:pf_Dhv-gv}
|\Delta_h v(x) - \nabla v(x)\cdot h\1_{\{|h|\le 1\}}| \le 
\begin{cases}
2\|v\|_{C^{\tilde{\sigma}}} \tilde{\sigma}(|h|\wedge 1), \quad &\text{ if } I_{\vphi\psi} \subset(1,2),\\
2(\|v\|_{C^0} + \|D^2v\|_{C^0})(|h|^2\wedge 1), \quad &\text{ if } I_{\vphi\psi}\subset(2,3),
\end{cases}
\end{align}
where $\tilde{\sigma}$ is a function on $(0,1]$ defined by 
$$ \tilde{\sigma}(t) = t^{(M_\vphi+m_\vphi+v_\psi)/2}.$$
Then we have from \eqref{l:aux_int} and \autoref{p:ineq_norms} that
\begin{align*}
 \|\cB v\|_{C^0} \le c_1 \|v\|_{C^0} + (\eps/5) \|v\|_{C^{\vphi\psi}}
 \end{align*}
for some constant $c_1 = c_1(r,\eps)>0$. 

For the $\psi$-H\"older seminorm of $\wt{\cB} v$, firstly we can use the inequality \eqref{e:pf_Bv_semi1} without any change. We can also use \eqref{e:pf_Bv_semi2} as it is, because $\nabla v(x+k) = \nabla v(x) = 0$ for $|k|\le r/2$ and $x\notin B(x_0,3r)$. Now suppose $x\in B(x_0,3r)$. This implies $|b(x+k,h)| \le c_2\psi(r)$ for any $h \in \R^d$ and $|k| \le r/2$. We have
\begin{align*}
\Delta_k(\wt{\cB}v) (x) &= \int_{\R^d} (\Delta_k\Delta_h v(x) - \Delta_k(\nabla v)(x)\cdot h \1_{\{|h|\le 1\}}) b(x+k,h) \frac{\d h }{|h|^d\vphi(|h|)}  \\
&\quad + \int_{\R^d} (\Delta_h v(x) - \nabla v(x)\cdot h \1_{\{|h|\le 1\}}) \Delta_k(b(\cdot,h))(x) \frac{\d h }{|h|^d\vphi(|h|)}\\
& = I_6 + I_7.
\end{align*}
If $I_{\vphi\psi} \subset(1,2)$, then we have
\begin{align}\label{e:pf_DkDhv-gv}
|\Delta_k\Delta_h v(x) - \Delta_k(\nabla v)(x)\cdot h \1_{\{|h|\le 1\}}| \le 2\cdot
\begin{cases}
\|v\|_{C^{\vphi\psi}} (\vphi\psi)(|h|), \quad & \text{ if } |h|\le |k|,\\
\|v\|_{C^{\vphi\psi}} |h|\frac{(\vphi\psi)(|k|)}{|k|}, \quad & \text{ if } |k|< |h| \le 1, \\
\|v\|_{C^\psi} \psi(|k|),\quad & \text{ if } |h|>1.
\end{cases}
\end{align}
Taking the integration, we get
\begin{align*}
|I_6| &\le \int_{|h|\le |k|} 2\|v\|_{C^{\vphi\psi}} (\vphi\psi)(|h|)\frac{c_2\psi(r)}{|h|^d\vphi(|h|)} \d h \\
& \quad + \int_{|k|< |h|\le 1} 2\|v\|_{C^{\vphi\psi}} \frac{|h|(\vphi\psi)(|k|)}{|k|}\cdot \frac{c_2\psi(r)}{|h|^d\vphi(|h|)} \d h \\
&\quad + \int_{|h|>1 } 2\|v\|_{C^\psi} \psi(|k|) \frac{c_2\psi(r)}{|h|^d\vphi(|h|)} \d h \\
&\le c_3 \psi(r)(\|v\|_{C^\psi}+\|v\|_{C^{\vphi\psi}}) \psi(|k|).
\end{align*}
For the case $1\in I_\psi$, the additional assumption that $a(x,h)=a(x,-h)$ for $x,h\in\R^d$ allows for the integrand over the region $\{h\in \R^d: |k|< |h|\le 1\}$ to be reduced to
$$|\Delta_k\Delta_h v(x)|\le 2\|v\|_{C^{\vphi\psi}} \frac{|k|(\vphi\psi)(|h|)}{|h|}.$$
If $I_{\vphi\psi} \subset(2,3)$, then we just replace \eqref{e:pf_DkDhv-gv} with 
\begin{align*}
|\Delta_k\Delta_h v(x) - \Delta_k(\nabla v)(x)\cdot h \1_{\{|h|\le 1\}}| \le 2\cdot
\begin{cases}
\|v\|_{C^{\vphi\psi}} \frac{|h|^2(\vphi\psi)(|k|)}{|k|^2}, \quad & \text{ if } |h|\le |k|,\\
\|v\|_{C^{\vphi\psi}} \frac{|k|(\vphi\psi)(|h|)}{|h|}, \quad & \text{ if } |k|< |h| \le 1, \\
\|v\|_{C^\psi} \psi(|k|),\quad & \text{ if } |h|>1.
\end{cases}
\end{align*}
Hence applying \autoref{p:ineq_norms} to $\|\cdot \|_{C^\psi}$ and $\|\cdot \|_{C^{\vphi\psi}}$ implies
\begin{align*}
|I_6| \le c_4\psi(r)(\|v\|_{C^0} + \|v\|_{C^{\vphi\psi}}) \psi(|k|),
\end{align*}
whenever $I_{\vphi\psi}\subset(1,2)\cup(2,3)$.

For $I_7$, we use $|\Delta_k (b(\cdot, h))(x)| \le 2\Lambda_3 \psi(|k|)$ and \eqref{e:pf_Dhv-gv} to get
\begin{align*}
|I_7| \le (c_5 \|v\|_{C^0} + (\eps/5)\|v\|_{C^{\vphi\psi}}) \psi(|k|),
\end{align*}
for some constant $c_5 = c_5(r,\eps)>0$. By summing up the above result and choosing $r$ such that $c_4\psi(r) \le \eps/5$, we get the result.
\qed

\begin{remark}
In the proof of \autoref{l:Bv} and \autoref{l:BBv}, the fact that $b(x,h)=a(x,h)-a(x_0,h)$ is only used to estimate $I_4$ and $I_6$. It allows us to find $r$ depending  on $\eps$. If we fix $\eps$ as a number, then we can replace $b(x,h)$ with $a(x,h)$ and obtain 
\begin{align}\label{e:Leta}
\| \cL \eta\|_{C^{\psi}} \le C.
\end{align}
for a constant $C$ depending on $r$.
\end{remark}

We are finally ready to prove our main result. 

\pf [Proof of \autoref{theo:main}] By \autoref{p:u_eta} it is enough to 
show that for any $\eps >0$ there exist positive constants $r$ and $c_1$ such that for all 
$x_0\in\bR^d$ 
\begin{align}\label{e:proof-auxiliary}
\|u\eta_{r,x_0} \|_{C^{\vphi\psi}} \le c_1\|f\|_{C^\psi} + c_1 \|u\|_{C^0} + 
\eps 
\|u\|_{C^{\vphi\psi}} \,.
\end{align}

First we consider the case $M_\vphi<1$. We define a freezing operator
\begin{align*}
\cL_{0}u(x) = \int_{\R^d} \Delta_hu(x) \frac{a(x_0,h)}{|h|^d\vphi(|h|)} \d h\,,
\end{align*}
and  $\cB = \cL-\cL_0$. 
Let $v(x)= u(x)\eta_{r,x_0}(x)$. As we mentioned at the beginning of this 
section we write $\eta$ instead of $\eta_{r,x_0}$. Observe that the identity
\begin{align*}
\Delta_h(u\eta)(x) = \eta(x)\Delta_h u(x) + u(x)\Delta_h\eta(x) +  \Delta_h 
u(x)\Delta_h\eta(x),
\end{align*}
yields
\begin{align*}
\cL v(x) = \eta(x) \cL u(x) +u(x) \cL \eta(x) +  H(x),
\end{align*}
where $H$ is defined by
$$ H(x) = \int_{\bR^d} 
\Delta_hu(x)\Delta_h\eta(x)\frac{a(x,h)}{|h|^d\vphi(|h|)}\d h.$$
Then we have
\begin{align*}
\cL_0 v(x) &= \eta(x)f(x) + u(x)\cL\eta(x) +  H(x) - \cB v(x).
\end{align*}
\autoref{theo:mainL_0-intro} now 
implies that with some constant $c_2 \geq 1$ 
\begin{align*}
\| v\|_{C^{\vphi\psi}} \le c_2( \|\eta f + u \cL \eta + H -\cB v  \|_{C^\psi} + \|v\|_{C^0}).
\end{align*}
Choose $r>0$ from \autoref{l:Bv} such that 
\begin{align*}
\|\cB v\|_{C^\psi} \le c_3\|v\|_{C^0} + (2c_2)^{-1} \|v\|_{C^{\vphi\psi}}.
\end{align*}
It is obvious that $\|\eta f\|_{C^{\psi}} \le c_4\|f\|_{C^\psi}$. It follows from \eqref{e:Leta} and \autoref{p:ineq_norms} that 
\begin{align*}
\|u \cL \eta\|_{C^\psi} \le c_6\|u\|_{C^0} + (4c_2)^{-1}\eps\|u\|_{C^{\vphi\psi}}.
\end{align*}
Finally, \eqref{l:H} implies that for a given $\eps>0$, there exists a constant $c_7=c_7(r,\eps)>0$ such that
\begin{align*}
\|H\|_{C^\psi} \le c_7 \|u\|_{C^0} + (4c_2)^{-1}\eps\|u\|_{C^{\vphi\psi}}.
\end{align*}
Hence \eqref{e:proof-auxiliary} holds true.
For the case $M_\vphi \ge 1$, we denote a freezing operator by 
\begin{align*}
\wt{\cL}_0 u (x) = \int_{\R^d} (\Delta_h u(x) -\nabla u(x)\cdot h \1_{\{|h|\le 1\}}) \frac{a(x_0,h)}{|h|^d\vphi(|h|)} \d h,
\end{align*}
and $\wt{\cB} = \cL - \wt{\cL}_0$. Then we have
\begin{align*}
\wt{\cL}_0 v(x) &= \eta(x)f(x) + u(x)\cL\eta(x) +  H(x) - \wt{\cB} v(x).
\end{align*}
Using \autoref{l:BBv} instead of \autoref{l:Bv}, we get the result from the same argument above. \qed

%%%%%%%%%%%%%%%%%%%%%%%%%%%%%%%%%%%%%%%%%%%%%%%%%%%%%%%%%%%%%%%%%%%%%%%%%%%%
%%%%%%%%%%%%%%%%%%%%%%%%%%%%%%%%%%%%%%%%%%%%%%%%%%%%%%%%%%%%%%%%%%%%%%%%%%%%
\section{Continuity of $\cL$}\label{sec:maprop}
%%%%%%%%%%%%%%%%%%%%%%%%%%%%%%%%%%%%%%%%%%%%%%%%%%%%%%%%%%%%%%%%%%%%%%%%%%%%
%%%%%%%%%%%%%%%%%%%%%%%%%%%%%%%%%%%%%%%%%%%%%%%%%%%%%%%%%%%%%%%%%%%%%%%%%%%%

In this section we provide the proof of \autoref{theo:maprop}. 

\pf [Proof of \autoref{theo:maprop}]  We first consider the case $M_\vphi <1$ that $\cL$ is defined by
$$ \cL u(x) = \int_{\R^d} (u(x+h)-u(x)) \frac{a(x,h)}{|h|^d\vphi(|h|)} \d h .$$
We observe that 
\begin{align*}
|\Delta_h u(x) | \le 2\|u\|_{C^{\vphi\psi}} 
\begin{cases}
(\vphi \psi)(|h|\wedge 1) ,\quad & I_{\vphi\psi}\subset (0,1),\\
|h|\wedge 1 , \quad & I_{\vphi\psi} \subset(1,2).
\end{cases}
\end{align*}
It follows from \autoref{l:aux_int} that $|\cL u(x)| \le c_1 \|u\|_{C^{\vphi\psi}}$ for all $x\in \R^d$. For the seminorm of $\cL u$, we know that 
\begin{align*}
\Delta_k (\cL u)(x)  = \int_{\R^d} \Delta_k\Delta_h u(x) \frac{a(x+k,h)}{|h|^d\vphi(|h|)} \d h + \int_{\R^d} \Delta_h u(x) \frac{\Delta_k (a(\cdot, h))(x)}{|h|^d\vphi(|h|)} \d h. 
\end{align*}
For the case $I_{\vphi\psi} \subset (0,1)$, we have
\begin{align*}
|\Delta_k\Delta_h u(x) | \le 2\|u\|_{C^{\vphi\psi}}(\vphi\psi)(|h|\wedge |k|),
\end{align*}
and, if $I_{\vphi\psi} \subset (1,2)$, then
\begin{align*}
|\Delta_k\Delta_h u(x)| \le 4\|u\|_{C^{\vphi\psi}} \cdot
\begin{cases}
\frac{(\vphi\psi)(|k|)}{|k|} |h| ,\quad & |h|\le |k|,\\
\frac{(\vphi\psi)(|h|)}{|h|} |k|, \quad & |k|<|h|\le 1,\\
\psi(|k|), \quad & |h|>1.
\end{cases}
\end{align*}
We used \autoref{p:ineq_norms} with the assumption in the last case. By 
integrating the above terms with respect to $h$ and using 
assumption \eqref{eq:assum_axh-cont} and \autoref{l:aux_int}, we get
$$ |\Delta_k(\cL u)(x)| \le c_2\|u\|_{C^{\vphi\psi}} \psi(|k|)$$
for every $x\in \R^d$ and $|k|\le 1$, which implies $[\cL u]_{C^{0;\psi}} \le c_2\|u\|_{C^{\vphi\psi}}$.

For the case $M_\vphi \ge 1$, the operator that we consider is given by 
\begin{align*}
\cL u(x) = \int_{\R^d} (u(x+h)-u(x)-\nabla u(x)\cdot h \1_{\{|h|\le 1\}}) \frac{a(x,h)}{|h|^d\vphi(|h|)} \d h. 
\end{align*}
When $1\notin I_\vphi$, similarly to the above case, we only need to observe that
\begin{align*}
|\Delta_h u(x)- \nabla u(x) \cdot h\1_{\{|h|\le 1\}} | \le 2\|u\|_{C^{\vphi\psi}} \cdot
\begin{cases}
(\vphi \psi)(|h|\wedge 1) ,\quad & I_{\vphi\psi}\subset (1,2),\\
(|h|\wedge 1)^2 , \quad & I_{\vphi\psi} \subset(2,3),
\end{cases}
\end{align*}
and if $I_{\vphi\psi} \subset(1,2)$,
\begin{align}\label{e:proof-aux_DkDhu}
|\Delta_k\Delta_h u(x) - \Delta_k \nabla u(x) \cdot h\1_{\{|h|\le 1\}} | \le c_3\|u\|_{C^{\vphi\psi}}\cdot
\begin{cases}
\frac{(\vphi\psi)(|h|\wedge |k|)}{|h|\wedge |k|}|h| , \quad & |h|\le 1,\\
\psi(|k|), \quad & |h|> 1,
\end{cases}
\end{align}
and if $I_{\vphi\psi} \subset(2,3)$,
\begin{align*}
|\Delta_k\Delta_h u(x) - \Delta_k \nabla u(x)\cdot h \1_{\{|h|\le 1\}} | \le c_4\|u\|_{C^{\vphi\psi}}\cdot
\begin{cases}
\frac{(\vphi\psi)(|k|)}{|k|^2} |h|^2 ,\quad & |h|\le |k|,\\
\frac{(\vphi\psi)(|h|)}{|h|^2} |h||k|, \quad & |k|< |h|\le 1, \\
\psi(|k|),\quad & |h|>1.
\end{cases}
\end{align*}
When $1\in I_\vphi$, it is easily shown $I_{\vphi\psi}\subset(1,2)$, and  we can replace the indicator function $\1_{\{|h|\le 1\}}$  in \eqref{e:proof-aux_DkDhu} with $\1_{\{|h|\le |k|\}}$ from the symmetry of $h \mapsto a(x, h)$. Thus we replace \eqref{e:proof-aux_DkDhu} with 
\begin{align*}
|\Delta_k\Delta_h u(x) - \Delta_k \nabla u(x) \cdot h\1_{\{|h|\le 1\}} | \le c_3\|u\|_{C^{\vphi\psi}}\cdot
\begin{cases}
(\vphi\psi)(|h|) , \quad & |h|\le |k|,\\
\frac{(\vphi\psi)(|h|)}{|h|} |k|, \quad &|k|<|h|\le 1,\\
\psi(|k|), \quad & |h|> 1.
\end{cases}
\end{align*}
Calculations similar to the case $M_\vphi<1$ give rise to the result. \qed

% \nocite*

\bibliographystyle{amsalpha}
\bibliography{literature_bae-kassmann}

\def\cprime{$'$} \def\cprime{$'$}
\providecommand{\bysame}{\leavevmode\hbox to3em{\hrulefill}\thinspace}
\providecommand{\MR}{\relax\ifhmode\unskip\space\fi MR }
\providecommand{\MRhref}[2]{%
  \href{http://www.ams.org/mathscinet-getitem?mr=#1}{#2}
}
\providecommand{\href}[2]{#2}
\begin{thebibliography}{KKK13}

\bibitem[Bas09]{Bas09}
R.~F. Bass, \emph{Regularity results for stable-like operators}, J. Funct.
  Anal. \textbf{257} (2009), no.~8, 2693--2722. \MR{2555009 (2010k:47099)}

\bibitem[BF75]{BeFo75}
C.~Berg and G.~Forst, \emph{Potential theory on locally compact abelian
  groups}, Springer-Verlag, New York-Heidelberg, 1975, Ergebnisse der
  Mathematik und ihrer Grenzgebiete, Band 87. \MR{0481057 (58 \#1204)}

\bibitem[BGT89]{BGT89}
N.~H. Bingham, C.~M. Goldie, and J.~L. Teugels, \emph{Regular variation},
  Encyclopedia of Mathematics and its Applications, vol.~27, Cambridge
  University Press, Cambridge, 1989. \MR{1015093 (90i:26003)}

\bibitem[BS56]{BaSt56}
N.~K. Bari and S.~B. Ste{\v{c}}kin, \emph{Best approximations and differential
  properties of two conjugate functions}, Trudy Moskov. Mat. Ob\v s\v c.
  \textbf{5} (1956), 483--522. \MR{0080797 (18,303e)}

\bibitem[DK13]{DoKi13}
H.~Dong and D.~Kim, \emph{Schauder estimates for a class of non-local elliptic
  equations}, Discrete Contin. Dyn. Syst. \textbf{33} (2013), no.~6,
  2319--2347. \MR{3007688}

\bibitem[FL06]{FaLe06}
W.~Farkas and H.-G. Leopold, \emph{Characterisations of function spaces of
  generalised smoothness}, Ann. Mat. Pura Appl. (4) \textbf{185} (2006), no.~1,
  1--62. \MR{2179581 (2008j:46024)}

\bibitem[Gol72]{Gol72}
M.~L. Gol'dman, \emph{On the inclusion of generalized h{\"o}lder classes},
  Mathematical notes of the Academy of Sciences of the USSR \textbf{12} (1972),
  no.~3, 626--631.

\bibitem[GT83]{GiTr83}
D.~Gilbarg and N.~S. Trudinger, \emph{Elliptic partial differential equations
  of second order}, second ed., Grundlehren der Mathematischen Wissenschaften
  [Fundamental Principles of Mathematical Sciences], vol. 224, Springer-Verlag,
  Berlin, 1983. \MR{737190 (86c:35035)}

\bibitem[JX14]{JiXi14}
T.~Jin and J.~Xiong, \emph{Schauder estimates for nonlocal fully nonlinear
  equations}, http://arxiv.org/pdf/1405.0758v3, 2014.

\bibitem[KK15]{KiKi15}
I.~Kim and K.-H. Kim, \emph{A {H}\"{o}lder regularity theory for a class of
  non-local elliptic equations related to subordinate {B}rownian motions}, Pot.
  Anal. (2015), DOI 10.1007/s11118-015-9490-5.

\bibitem[KKK13]{KKK13}
I.~Kim, K.-H. Kim, and P.~Kim, \emph{Parabolic {L}ittlewood-{P}aley inequality
  for {$\phi(-\Delta)$}-type operators and applications to stochastic
  integro-differential equations}, Adv. Math. \textbf{249} (2013), 161--203.
  \MR{3116570}

\bibitem[KL87]{KiLi87}
G.~A. Kaljabin and P.~I. Lizorkin, \emph{Spaces of functions of generalized
  smoothness}, Math. Nachr. \textbf{133} (1987), 7--32. \MR{912417 (89b:46046)}

\bibitem[KN12]{KrNi12}
D.~Kreit and S.~Nicolay, \emph{Some characterizations of generalized {H}\"older
  spaces}, Math. Nachr. \textbf{285} (2012), no.~17-18, 2157--2172.
  \MR{3002607}

\bibitem[ROS14]{RoSe14}
X.~Ros-Oton and J.~Serra, \emph{Regularity theory for general stable
  operators}, http://arxiv.org/pdf/1412.3892v1, 2014.

\bibitem[Ul{\cprime}68]{Ulj68}
P.~L. Ul{\cprime}janov, \emph{The embedding of certain classes
  {$H_{p}{}^{\omega }$} of functions}, Izv. Akad. Nauk SSSR Ser. Mat.
  \textbf{32} (1968), 649--686. \MR{0231194 (37 \#6749)}

\end{thebibliography}

\end{document}